\documentclass[12pt]{amsart}
\usepackage{amssymb,latexsym}

\numberwithin{equation}{section}

\newtheorem{theorem}{Theorem}[section]
\newtheorem{proposition}[theorem]{Proposition}
\newtheorem{conjecture}[theorem]{Conjecture}
\newtheorem{corollary}[theorem]{Corollary}
\newtheorem{lemma}[theorem]{Lemma}

\theoremstyle{definition}

\newtheorem{remark}[theorem]{Remark}
\newtheorem{example}[theorem]{Example}
\newtheorem{definition}[theorem]{Definition}

\def\proof{\smallskip\noindent {\bf Proof. }}
\def\endproof{\hfill$\square$\medskip}

\newcommand{\doublesubscript}[3]{
\displaystyle\mathop{\displaystyle #1_{#2}}_{#3} }

\def\AA{\mathcal{A}}

\def\FF{\mathcal{F}}
\def\ZZ{\mathbb{Z}}
\def\QQ{\mathbb{Q}}
\def\CC{\mathbb{C}}
\def\RR{\mathbb{R}}
\def\LL{\mathcal{L}}

\def\OO{\mathcal{O}}
\def\gg{\mathfrak{g}}

\def\hh{\mathfrak{h}}
\def\ii{\mathbf{i}}

\def\cc{\mathbf{c}}
\def\xx{\mathbf{x}}
\def\TT{\mathcal{T}}
\def\ex{\mathbf{ex}}

\def\ZZ{\mathbb{Z}}
\def\PP{\mathbb{P}}

\newcommand{\mat}[4]{\left(\!\!\begin{array}{cc}
#1 & #2 \\ #3 & #4 \\
\end{array}\!\!\right)}

\begin{document}

\title{Quantum cluster algebras}

\author{Arkady Berenstein}
\address{\noindent Department of Mathematics, University of Oregon,
Eugene, OR 97403, USA} \email{arkadiy@math.uoregon.edu}

\author{Andrei Zelevinsky}
\address{\noindent Department of Mathematics, Northeastern University,
  Boston, MA 02115, USA}
\email{andrei@neu.edu}

\date{May 7, 2004; revised June 28, 2004}

\thanks{Research supported in part
by NSF (DMS) grants \# 0102382 (A.B.) and 0200299~(A.Z.).}


\maketitle

\makeatletter
\renewcommand{\@evenhead}{\tiny \thepage \hfill  A.~BERENSTEIN and  A.~ZELEVINSKY \hfill}

\renewcommand{\@oddhead}{\tiny \hfill QUANTUM CLUSTER ALGEBRAS
 \hfill \thepage}
\makeatother

\tableofcontents


\section{Introduction}
\label{sec:introduction}

Cluster algebras were introduced by S.~Fomin and A.~Zelevinsky in
\cite{fz-clust1}; their study continued in \cite{fz-clust2,bfz-clust3}.
This is a family of commutative rings designed to serve as an
algebraic framework for the theory of total positivity and
canonical bases in semisimple groups and their quantum analogs.
In this paper we introduce and study quantum deformations of cluster algebras.

Our immediate motivation for introducing quantum cluster algebras
is to prepare the ground for a general notion of the canonical
basis in a cluster algebra.
Remarkably, cluster algebras and their quantizations appear to be
relevant for the study of (higher) Teichmuller theory initiated
in \cite{gsv,gsv2,FG1,FG2}.
Our approach to quantization has much in common with the one in \cite{FG1,FG2},
but we develop it more systematically.
In particular, we show that practically all the structural results on cluster
algebras obtained in \cite{fz-clust1,fz-clust2,bfz-clust3} extend
to the quantum setting.
This includes the Laurent phenomenon  \cite{fz-clust1,fz-laurent,bfz-clust3}
and the classification of cluster algebras of finite type \cite{fz-clust2}.

Our approach to quantum cluster algebras can be described as follows.
Recall that a \emph{cluster algebra}~$\AA$ is a certain commutative ring generated
by a (possibly infinite) set of generators called \emph{cluster variables}
inside an ambient field~$\FF$ isomorphic to the field of rational
functions in~$m$ independent variables over~$\QQ$.
The set of cluster variables is the union of some distinguished
transcendence bases of~$\FF$ called (extended) \emph{clusters}.
The clusters are not given from the outset but are obtained from
an initial cluster via an iterative process of \emph{mutations}
which follows a set of canonical rules.
According to these rules, every cluster $\{x_1, \dots, x_m\}$
is surrounded by~$n$ adjacent clusters (for some $n \leq m$ called the \emph{rank}
of~$\AA$) of the form $\{x_1, \dots, x_m\} - \{x_k\} \cup \{x'_k\}$,
where~$k$ runs over a given~$n$-element subset of
\emph{exchangeable} indices, and $x'_k \in \FF$ is related to $x_k$
by the \emph{exchange relation} (see \eqref{eq:exchange-rel-xx} below).
The cluster algebra structure is completely determined by an $m
\times n$ integer matrix~$\tilde B$ that encodes all the exchange relations.
(The precise definitions of all these notions are given in
Section~\ref{sec:geom-type} below.)
Now the quantum deformation of~$\AA$ is a $\QQ(q)$-algebra obtained by making each
cluster into a \emph{quasi-commuting} family $\{X_1, \dots, X_m\}$;
this means that $X_i X_j = q^{\lambda_{ij}} X_j X_i$ for a
skew-symmetric integer $m \times m$ matrix $\Lambda = (\lambda_{ij})$.
In doing so, we have to modify the mutation process and the
exchange relations so that all the adjacent quantum clusters will
also be quasi-commuting.
This imposes the \emph{compatibility} relation between the
quasi-commutation matrix~$\Lambda$ and the exchange matrix~$\tilde B$
(Definition~\ref{def:compatible-triple} below).
In what follows, we develop a formalism that allows us to show
that any compatible matrix pair $(\Lambda, \tilde B)$
gives rise to a well defined quantum cluster algebra.

The paper is organized as follows.
In Section~\ref{sec:geom-type} we
present necessary definitions and facts from the theory of
cluster algebras in the form suitable for our current purposes.
In Section~\ref{sec:strongly-skew-symm}, we introduce compatible
matrix pairs $(\Lambda, \tilde B)$ and their mutations.

Section~\ref{sec:qdeformations} plays the central part in this paper.
It introduces the main concepts needed for the definition of
quantum cluster algebras (Definition~\ref{def:quantum-cluster-algebra}):
\emph{based quantum tori} (Definition~\ref{def:quantum-torus}) and their
skew-fields of fractions, \emph{toric frames}
(Definition~\ref{def:toric-frame}), \emph{quantum seeds}
(Definition~\ref{def:quantum-seed}) and their mutations
(Definition~\ref{def:quantum-seed-mutation}).

Section~\ref{sec:upper-bound} establishes the quantum version of
the Laurent phenomenon (Corollary~\ref{cor:A-subset-U}): any
cluster variable is a Laurent polynomial in the elements of any
given cluster.
The proof closely follows the argument in \cite{bfz-clust3}
with necessary modifications.
It is based on the important concept of an \emph{upper cluster
algebra} and the fact that it is invariant under mutations
(Theorem~\ref{th:upper-bound-invariance}).

In Section~\ref{sec:exchange-graphs}, we show that the
\emph{exchange graph} of a quantum cluster algebra remains unchanged
in the ``classical limit" $q = 1$
(Theorem~\ref{th:exchange-graphs}).
(Recall that the vertices of the exchange graph correspond to
(quantum) seeds, and the edges correspond to mutations.)
An important consequence of Theorem~\ref{th:exchange-graphs} is
that the classification of cluster algebras
of finite type achieved in \cite{fz-clust2} applies verbatim to
quantum cluster algebras.

An important ingredient of the proof of Theorem~\ref{th:exchange-graphs}
is the \emph{bar-involution} on the quantum cluster algebra which
is modeled on the Kazhdan-Lusztig involution, or the one used later by Lusztig in
his definition of the canonical basis.
We conclude Section~\ref{sec:exchange-graphs} by including the
bar-involution into a family of \emph{twisted bar-involutions}
(Proposition~\ref{pr:twisted-bar-involution}).
This construction is  motivated by our hope that this family of
involutions will find applications to the future theory of
canonical bases in (quantum) cluster algebras.

Section~\ref{sec:lower} extends to the quantum setting another
important result from \cite{bfz-clust3}: a sufficient condition
(``acyclicity") guaranteeing that the cluster algebra coincides
with the upper one (Theorem~\ref{th:acyclic-closing-gaps}).
The proof in \cite{bfz-clust3} is elementary but rather involved;
we do not reproduce it here in the quantum setting, just indicate
necessary modifications.

Section~\ref{sec:Cartan-triples} presents our main source of
examples of quantum cluster algebras: those associated with double
Bruhat cells in semisimple groups.
The ordinary cluster algebra structure associated with these cells
was introduced and studied in \cite{bfz-clust3}.
The main result in Section~\ref{sec:Cartan-triples}
(Theorem~\ref{th:compatible-ii}) shows in particular that every matrix
$\tilde B$ associated as in \cite{bfz-clust3} with a double Bruhat
cell can be naturally included into a compatible matrix pair
$(\Lambda, \tilde B)$.
Not very surprisingly, the skew-symmetric matrix~$\Lambda$ that
appears here is the one describing the standard Poisson structure
in the double cell in question; this matrix was calculated in
\cite{kogzel,gsv}.
The statement and proof of Theorem~\ref{th:compatible-ii} are
purely combinatorial, i.e., do not use the geometry of double
cells; thus, without any additional difficulty, we state and prove it in greater
generality that allows us to produce a substantial class of
compatible matrix pairs associated with generalized Cartan
matrices.

The study of quantum double Bruhat cells continues in
Section~\ref{sec:quantum double Bruhat cells}.
(For the convenience of the reader, we collect necessary
preliminaries on quantum groups in Section~\ref{sec:quantum groups}.)
The goal is to relate the cluster algebra approach with that
developed by De Concini and Procesi in
\cite{dc-pr} (see also \cite{joseph,brown-goodearl}).
Our results here are just the first step in this direction; we merely prepare the
ground for a conjecture (Conjecture~\ref{con:birational isomorphism})
that every quantum double Bruhat cell is naturally isomorphic to
the upper cluster algebra associated with an appropriate matrix
pair from Theorem~\ref{th:compatible-ii}.
The classical case of this conjecture was proved in
\cite[Theorem~2.10]{bfz-clust3}.

For the convenience of the reader, some needed facts on Ore
localizations are collected with proofs in concluding
Section~\ref{sec:appendix}.

\section{Cluster algebras of geometric type}
\label{sec:geom-type}

We start by recalling the definition of (skew-symmetrizable) cluster
algebras of geometric type, in the form most convenient for our
current purposes.

Let $m$ and $n$ be two positive integers with $m \geq n$.
Let $\FF$ be the field of rational
functions over $\QQ$ in $m$ independent (commuting) variables.
The cluster algebra that we are going to introduce will be a
subring of the ambient field~$\FF$.
To define it, we need to introduce seeds and their mutations.

\begin{definition}
\label{def:seed}
A (skew-symmetrizable) \emph{seed} in $\FF$ is a pair
$(\tilde \xx, \tilde B)$, where
\begin{enumerate}

\item $\tilde \xx = \{x_1, \dots, x_m\}$
is a transcendence basis of~$\FF$, which generates~$\FF$.

\item $\tilde B$ is an $m \times n$ integer matrix with rows
labeled by $[1,m] = \{1, \dots, m\}$
and columns labeled by an $n$-element subset
$\ex \subset [1,m]$,
such that the $n \times n$ submatrix $B$ of $\tilde B$ with rows
labeled by $\ex$ is skew-symmetrizable, i.e.,
$DB$ is skew-symmetric for some diagonal $n \times n$ matrix $D$ with positive
diagonal entries.
\end{enumerate}
The seeds are defined up to a relabeling of elements of $\tilde \xx$
together with the corresponding relabeling of rows and columns of $\tilde B$.
\end{definition}

\begin{remark}
The last condition in~(1), namely that $\tilde \xx$
generates~$\FF$, was unfortunately omitted in
\cite{fz-clust2,bfz-clust3} although it was always meant to be there.
(We thank E.~B.~Vinberg for pointing this out to us.)
In what follows, we refer to the subsets satisfying (1) as
\emph{free generating sets} of~$\FF$.
\end{remark}

We denote $\xx = \{x_j: j \in \ex\} \subset \tilde \xx$, and
$\cc = \tilde \xx - \xx$.
We refer to the indices from $\ex$ as \emph{exchangeable indices},
to $\xx$ as the \emph{cluster} of a seed $(\tilde \xx, \tilde B)$, and to
$B$ as the \emph{principal part} of $\tilde B$ .

Following \cite[Definition~4.2]{fz-clust1}, we say that a real
$m \times n$ matrix $\tilde B'$ is obtained from $\tilde B$ by
\emph{matrix mutation} in direction~$k \in \ex$,
and write $\tilde B' = \mu_k (\tilde B)$
if the entries of $\tilde B'$ are given by
\begin{equation}
\label{eq:matrix-mutation}
b'_{ij} =
\begin{cases}
-b_{ij} & \text{if $i=k$ or $j=k$;} \\[.05in]
b_{ij} + \displaystyle\frac{|b_{ik}| b_{kj} +
b_{ik} |b_{kj}|}{2} & \text{otherwise.}
\end{cases}
\end{equation}
This operation has the following properties.

\begin{proposition}
\label{pr:matrix-mutation}
\begin{enumerate}
\item The principal part of $\tilde B'$ is equal to $\mu_k (B)$.

\item $\mu_k$ is involutive: $\mu_k (\tilde B') = \tilde B$.

\item If $B$ is integer and skew-symmetrizable then so is $\mu_k (B)$.

\item The rank of $\tilde B'$ is equal to the rank of $\tilde B$.
\end{enumerate}
\end{proposition}

\begin{proof}
Parts (1) and (2) are immediate from the definitions.
To see (3), notice that $\mu_k (B)$ has the same skew-symmetrizing
matrix $D$ (see \cite[Proposition~4.5]{fz-clust1}).
Finally, Part (4) is proven in \cite[Lemma~3.2]{bfz-clust3}.
\end{proof}

\begin{definition}
\label{def:seed-mutation}
Let $(\tilde \xx, \tilde B)$ be a seed in $\FF$. 
For any exchangeable index $k$, the \emph{seed mutation} in direction~$k$
transforms $(\tilde \xx, \tilde B)$ into a seed
$\mu_k(\tilde \xx, \tilde B)=(\tilde \xx', \tilde B')$, where
\begin{itemize}

\item
$\tilde \xx' = \tilde \xx - \{x_k\} \cup \{x'_k\}$, where $x'_k \in \FF$ is determined
by the \emph{exchange relation}
\begin{equation}
\label{eq:exchange-rel-xx}
x_k\,x'_k
= \prod_{\substack{i \in [1,m] \\ b_{ik}>0}} x_i^{b_{ik}}
+ \prod_{\substack{i \in [1,m] \\ b_{ik}<0}} x_i^{-b_{ik}}.
\end{equation}

\item
The matrix $\tilde B'$ is obtained from $\tilde B$ by the matrix
mutation in direction~$k$.

\end{itemize}
\end{definition}

Note that $(\tilde \xx', \tilde B')$ is indeed a seed, since $\tilde \xx'$ is
obviously a free generating set for $\FF$, 
and the principal part
of $\tilde B'$ is skew-symmetrizable by parts (1) and (3) of
Proposition~\ref{pr:matrix-mutation}.
As an easy consequence of part (2) of Proposition~\ref{pr:matrix-mutation},
the seed mutation is involutive, i.e.,
$\mu_k (\tilde \xx', \tilde B') = (\tilde \xx, \tilde B)$.
Therefore, the following relation on seeds is an equivalence
relation: we say that $(\tilde \xx, \tilde B)$ is mutation-equivalent to
$(\tilde \xx', \tilde B')$ and write
$(\tilde \xx, \tilde B) \sim (\tilde \xx', \tilde B')$ if $(\tilde \xx', \tilde B')$
can be obtained from $(\tilde \xx, \tilde B)$ by a sequence of seed mutations.
Note that all seeds $(\tilde \xx', \tilde B')$ mutation-equivalent to
a given seed $(\tilde \xx, \tilde B)$ share the same set
$\cc = \tilde \xx' - \xx'$.
Let $\ZZ[\cc^{\pm 1}] \subset \FF$ be the ring of integer
Laurent polynomials in the elements of $\cc$.

Now everything is in place for defining cluster algebras.

\begin{definition}
\label{def:cluster-algebra}
Let $\mathcal{S}$ be a mutation-equivalence class of seeds in $\FF$.
The cluster algebra $\AA(\mathcal{S})$ associated with $\mathcal{S}$
is the $\ZZ[\cc^{\pm 1}]$-subalgebra of the ambient field $\FF$,
generated by the union of clusters of all seeds in $\mathcal{S}$.
\end{definition}

Since $\mathcal{S}$ is uniquely determined by each of the seeds
$(\tilde \xx, \tilde B)$ in it, we sometimes
denote $\AA(\mathcal{S})$ as $\AA(\tilde \xx,\tilde B)$, or even simply $\AA(\tilde B)$,
because $\tilde B$ determines this algebra uniquely up to an
automorphism of the ambient field~$\FF$.

\section{Compatible pairs}
\label{sec:strongly-skew-symm}

\begin{definition}
\label{def:compatible-triple}
Let $\tilde B$ be an $m \times n$ integer matrix with rows labeled
by $[1,m]$ and columns labeled by an $n$-element subset
$\ex \subset [1,m]$.
Let $\Lambda$ 
be a skew-symmetric
$m \times m$ integer matrix with rows and columns
labeled by~$[1,m]$.
We say that a pair $(\Lambda, \tilde B)$ is
\emph{compatible} if, 
for every $j \in \ex$ and $i \in [1,m]$,
we have
$$\sum_{k = 1}^m b_{kj} \lambda_{ki} = \delta_{ij} d_j$$
for some positive integers $d_j \,\, (j \in \ex)$.
In other words, the $n \times m$ matrix
$\tilde D = \tilde B^T \Lambda$ consists of the two blocks:
the $\ex \times \ex$ diagonal matrix $D$ with positive integer
diagonal entries $d_j$, and the $\ex \times ([1,m] - \ex)$ zero block.
\end{definition}

A large class of compatible pairs is constructed in
Section~\ref{sec:cartan-data} below.
Here is one specific example of a pair from this class.

\begin{example}
\label{ex:sl3-quantized-matrices}
Let~$\tilde B$ be an $8 \times 4$ matrix given by
$$\tilde B =
  \begin{pmatrix}
-1 &  0 &  0 &  0\\
1 &  -1 &  0 &  0\\
0 & 1 &  -1 &  0\\
-1 &  0 & 1 &  -1\\
1 &  -1 &  0 & 1\\
0 & 1 &  -1 &  0\\
0 &  -1 &  0 & 1\\
0 &  0 &  0 &  -1
  \end{pmatrix},
$$ where the columns are indexed by the set $\ex = \{3,4,5,6\}$
(note that the $4 \times 4$ submatrix of~$\tilde B$ on the rows
$\{3,4,5,6\}$ is skew-symmetric).
(This matrix describes the
cluster algebra structure in the coordinate ring of $SL_3$
localized at the four minors $\Delta_{1,3}$, $\Delta_{3,1}$,
$\Delta_{12,23}$, and $\Delta_{23,12}$; it is obtained from the
one in \cite[Figure~2]{bfz-clust3} by interchanging the first two
rows and changing the sign of all entries.)
Let us define a skew-symmetric $8 \times 8$ matrix $\Lambda$ by
$$\Lambda =
  \begin{pmatrix}
0 &  0 &  -1 &  -1 & -1 & 0 & 0 & 0\\
0 & 0 & 0 & -1 &  -1 & -1&  0 &  0\\
1 & 0 & 0 & 0 & -1 & 0 &  1 &  0\\
1 & 1 &  0 & 0 & 0 & 0 & 1 &  1\\
1 &  1 & 1 & 0 & 0 & 1 & 1 & 1\\
0 & 1 & 0 & 0 & -1 &  0 & 0 & 1\\
0 & 0 &  -1 & -1 & -1 & 0 & 0 & 0\\
0 &  0 &  0 &  -1 & -1 & -1 & 0 & 0
  \end{pmatrix}.
$$
A direct check shows that the pair $(\Lambda,\tilde B)$ is
compatible: the product $\tilde D = \tilde B^T \Lambda$ is equal to
$$\begin{pmatrix}
0 &  0 &  2 &  0 & 0 & 0 & 0 & 0\\
0 & 0 & 0 & 2 &  0 & 0&  0 &  0\\
0 & 0 & 0 & 0 & 2 & 0 &  0 &  0\\
0 & 0 &  0 & 0 & 0 & 2 & 0 &  0
  \end{pmatrix}.
$$
\end{example}

\begin{proposition}
\label{pr:stronger}
If a pair $(\Lambda, \tilde B)$ is compatible
then $\tilde B$ has full rank~$n$, and its
principal part $B$ is skew-symmetrizable.
\end{proposition}

\begin{proof}
By the definition, the $n \times n$ submatrix of $\tilde B^T \Lambda$
with rows and columns labeled by $\ex$ is
the  diagonal matrix $D$ with positive diagonal entries $d_j$.
This implies at once that ${\rm rk} (\tilde B) = n$.
To show that $B$ is skew-symmetrizable, note
that $DB = \tilde B^T \Lambda \tilde B$ is skew-symmetric.
\end{proof}

We will extend matrix mutations to those of compatible pairs.
Fix an index $k \in \ex$ and a sign $\varepsilon \in \{\pm 1\}$.
As shown in \cite[(3.2)]{bfz-clust3}, the matrix $\tilde B' = \mu_k(\tilde B)$
can be written as
\begin{equation}
\label{eq:mutation-product}
\tilde B' =
E_\varepsilon\,\tilde B\, F_\varepsilon
\, ,
\end{equation}
where
\begin{itemize}
\item
$E_\varepsilon$
is the $m\times m$ matrix
with entries
\begin{equation}
\label{eq:E-entries}
e_{ij} =
\begin{cases}
\delta_{ij} & \text{if $j \neq k$;} \\[.05in]
- 1 & \text{if $i = j = k$;} \\[.05in]
\max(0, -\varepsilon b_{ik})& \text{if $i \neq j = k$.}
\end{cases}
\end{equation}

\item
$F_\varepsilon$
is the $n\times n$ matrix with rows and columns labeled by $\ex$,
and entries given by
\begin{equation}
\label{eq:F-entries}
f_{ij} =
\begin{cases}
\delta_{ij} & \text{if $i \neq k$;} \\[.05in]
- 1 & \text{if $i = j = k$;} \\[.05in]
\max(0, \varepsilon b_{kj})& \text{if $i = k \neq j$.}
\end{cases}
\end{equation}
\end{itemize}

Now suppose that a pair $(\Lambda, \tilde B)$ is compatible.
We set
\begin{equation}
\label{eq:Lambda-mutated}
\Lambda' = E_\varepsilon^T \Lambda E_\varepsilon;
\end{equation}
thus, $\Lambda'$ is skew-symmetric.

\begin{proposition}
\label{pr:triple-mutation}
\begin{enumerate}
\item The pair $(\Lambda', \tilde B')$ is compatible.

\item $\Lambda'$ is independent of the choice of a sign~$\varepsilon$.
\end{enumerate}
\end{proposition}

\begin{proof}
To prove (1), we show that the pair $(\Lambda', \tilde B')$
satisfies Definition~\ref{def:compatible-triple}
with the same matrix~$\tilde D$.
We start with an easy observation that
\begin{equation}
\label{eq:EF-squares}
E_\varepsilon^2 = 1, \quad
F_\varepsilon^2 = 1.
\end{equation}
We also have
\begin{equation}
\label{eq:FTD}
 F_\varepsilon^T \tilde D = \tilde D E_\varepsilon;
 \end{equation}
indeed, one only has to check that
$$d_{i} \max(0, -\varepsilon b_{ik}) =
d_k \max(0, \varepsilon b_{ki})$$
for $i \in \ex - \{k\}$, which is true since,
by Proposition~\ref{pr:stronger},
$D$ is a skew-symmetrizing matrix for the principal part of $\tilde B$.
In view of \eqref{eq:EF-squares} and \eqref{eq:FTD}, we have
$$
(\tilde B')^T \Lambda' = F_\varepsilon^T \tilde D E_\varepsilon =
\tilde D,$$
finishing the proof.

\smallskip

(2) An easy calculation shows that the matrix entries of the product
$G = E_- E_+$ are given by
\begin{equation}
\label{eq:EE-matrix}
g_{ij} =
\begin{cases}
1 & \text{if $i=j$;} \\[.05in]
\delta_{jk} b_{ik}
& \text{if $i \neq j$.}
\end{cases}
\end{equation}
A direct check now shows that $G^T \Lambda G = \Lambda$.
(For instance, if $j \neq k$ then the $(k,j)$ entry of $G^T \Lambda G$
is equal to
$$\lambda_{kj} + \sum_{i \neq k} b_{ik} \lambda_{ij} = \lambda_{kj},$$
since the sum $\sum_{i \neq k} b_{ik} \lambda_{ij}$ is the $(k,j)$-entry of
$\tilde B^T \Lambda$ and so is equal to~$0$.)
We conclude that $E_+^T \Lambda E_+ = E_-^T \Lambda E_-$
as claimed.
\end{proof}

Proposition~\ref{pr:triple-mutation} justifies
the following important definition.

\begin{definition}
\label{def:triple-mutation}
Let $(\Lambda, \tilde B)$ be a compatible pair, and
$k \in \ex$.
We say that the compatible pair given by \eqref{eq:mutation-product}
and \eqref{eq:Lambda-mutated} is obtained from $(\Lambda, \tilde B)$
by the \emph{mutation} in direction~$k$, and write
$(\Lambda', \tilde B') = \mu_k (\Lambda, \tilde B)$.
\end{definition}

The following result extends part (2) of Proposition~\ref{pr:matrix-mutation}
to compatible pairs.

\begin{proposition}
\label{pr:triple-mutation-involutive}
The mutations of compatible pairs are involutive:
for any compatible pair $(\Lambda, \tilde B)$ and $k \in
\ex$, we have
$\mu_k (\mu_k (\Lambda, \tilde B)) = (\Lambda, \tilde B)$.
\end{proposition}

\proof
Let $\mu_k (\Lambda, \tilde B) = (\Lambda', \tilde B')$, and let
$E'_\varepsilon$ be given by \eqref{eq:E-entries} applied to
$\tilde B'$ instead of $\tilde B$.
By the first case in \eqref{eq:matrix-mutation}, the $k$th column
of $\tilde B'$ is the negative of the $k$th column of $\tilde B$.
It follows that
\begin{equation}
\label{E'-E}
E'_\varepsilon = E_{-\varepsilon}.
\end{equation}
In view of \eqref{eq:EF-squares}, we get
$$(E'_+)^T \Lambda' E'_+ = E_-^T \Lambda' E_- = \Lambda,$$
which proves the desired claim.
\endproof

\section{Quantum cluster algebras setup}
\label{sec:qdeformations}

\subsection{Based quantum torus and ambient skew-field}
Let $L$ be a lattice of rank $m$, with a 
skew-symmetric bilinear form $\Lambda : L\times L \to \ZZ$.
We also introduce a formal variable~$q$.
It will be convenient to work over the field of
rational functions $\QQ(q^{1/2})$ as a ground field.
Let $\ZZ[q^{\pm 1/2}] \subset \QQ(q^{1/2})$ denote the ring
of integer Laurent polynomials in the variable $q^{1/2}$.

\begin{definition}
\label{def:quantum-torus} The \emph{based quantum torus}
associated with $L$ is the $\ZZ[q^{\pm 1/2}]$-algebra
$\TT = \TT(\Lambda)$
with a distinguished $\ZZ[q^{\pm 1/2}]$-basis $\{X^e: e \in L\}$ and the
multiplication given by
\begin{equation}
\label{eq:multiplication-F}
X^e X^f = q^{\Lambda(e,f)/2} X^{e+f} \quad (e,f \in L) \, .
\end{equation}
\end{definition}

Thus, $\TT$ can be viewed as the group algebra
of $L$ over $\ZZ[q^{\pm 1/2}]$ twisted by a $2$-cocycle
$(e,f) \mapsto q^{\Lambda(e,f)/2}$.
It is easy to see that $\TT$ is associative: we have
\begin{equation}
\label{eq:multiplication-F-ass}
(X^e X^f) X^g = X^e (X^f X^g) =
q^{(\Lambda (e,f) + \Lambda (e,g) + \Lambda (f,g))/2} X^{e+f+g} \, .
\end{equation}
The basis elements satisfy the commutation relations
\begin{equation}
\label{eq:multiplication-F-com}
X^e X^f = q^{\Lambda (e,f)} X^f X^e \, .
\end{equation}
We also have
\begin{equation}
\label{eq:unit}
X^0=1, \quad (X^e)^{-1} = X^{-e} \quad (e \in L) \, .
\end{equation}

It is well-known (see Section~\ref{sec:appendix} below) that $\TT$ is an Ore
domain, i.e., is contained in its skew-field of fractions $\FF$.
Note that $\FF$ is a $\QQ(q^{1/2})$-algebra.
A quantum cluster algebra to be defined below will be a
$\ZZ[q^{\pm 1/2}]$-subalgebra of $\FF$.


\subsection{Some automorphisms of $\FF$}
\label{sec:auto}

Unless otherwise stated, by an \emph{automorphism} of $\FF$ we
will always mean a $\QQ(q^{1/2})$-algebra automorphism.
An important class of automorphisms of $\FF$ can be given as follows.
For a lattice point $b \in L - \ker (\Lambda)$,
let $d(b)$ denote the minimal positive value of $\Lambda(b,e)$
for $e \in L$.
We associate with $b$ the grading on $\TT$ such that every $X^e$
is homogeneous of degree
\begin{equation}
\label{eq:b-degree}
d_b (X^e) = d_b (e) = \Lambda(b,e)/d(b).
\end{equation}

\begin{proposition}
\label{pr:elementary-auto}
For every $b\in L - \ker (\Lambda)$,
and every sign $\varepsilon$, there is a unique automorphism
$\rho_{b,\varepsilon}$ of $\FF$ such that
\begin{equation}
\label{eq:rho}
\rho_{b,\varepsilon}(X^e) =
\begin{cases}
X^e & \text {if $\Lambda(b,e)=0$;}\\[.05in]
X^e+X^{e+\varepsilon b} & \text
{if $\Lambda(b,e)= -d(b)$.}
\end{cases}
\end{equation}
\end{proposition}

\begin{proof}
Since the elements $X^e$ that appear in \eqref{eq:rho}, together
with their inverses generate $\TT$ as a $\ZZ[q^{\pm 1/2}]$-algebra, the
uniqueness of $\rho_{b,\varepsilon}$ is clear.
To show the existence, we introduce some notation.
For every nonnegative integer $r$, we define an element
$P^r_{b,\varepsilon} \in \TT$ by
\begin{equation}
\label{eq:P-coeff}
P^r_{b,\varepsilon} = \prod_{p=1}^r (1 + q^{\varepsilon (2p-1)d(b)/2}
X^{\varepsilon b}).
\end{equation}
We extend the action of $\rho_{b,\varepsilon}$ given by
\eqref{eq:rho} to a $\ZZ[q^{\pm 1/2}]$-linear map
$\TT \to \FF$ such that, for every $e\in L$ with $|d_b (e)| = r$, we have
\begin{equation}
\label{eq:rho-extended}
\rho_{\varepsilon,b} (X^e) =
\begin{cases}
P^r_{b, \varepsilon} X^e   & \text{if $d_b (e) = -r$,} \\[.05in]
(P^r_{-b, -\varepsilon})^{-1} X^e  & \text{if $d_b (e) = r$.}
\end{cases}
\end{equation}
(it is easy to see that \eqref{eq:rho-extended} specializes to
\eqref{eq:rho} when $d_b (e) = 0$, or $d_b (e) = -1$; a more general
expression is given by \eqref{eq:rho-explicit} below).
One checks easily with the help of
\eqref{eq:multiplication-F-com} that this extended map is a $\ZZ[q^{\pm 1/2}]$-algebra
homomorphism $\TT \to \FF$, and so it extends to an algebra endomorphism of $\FF$.
The fact that this is an automorphism follows from the identity
$\rho_{-b, -\varepsilon} (\rho_{b,\varepsilon}(X^e)) = X^e$, which
is a direct consequence of \eqref{eq:rho-extended}.
\end{proof}

A direct check using \eqref{eq:rho-extended} shows that the
automorphisms $\rho_{b,\varepsilon}$ have the following
properties:
\begin{equation}
\label{eq:rho-mult}
\rho_{b,\varepsilon}^{-1} = \rho_{-b,-\varepsilon}, \quad
\rho_{b,-\varepsilon} = \tau_{b,\varepsilon} \circ \rho_{b,\varepsilon}  \ ,
\end{equation}
where $\tau_{b,\varepsilon}$ is an automorphism of $\FF$ acting by
$$\tau_{b,\varepsilon}(X^e) = X^{e - \varepsilon d_b(e) b} \quad (e \in L) \ .$$

In the first case in \eqref{eq:rho-extended}, i.e., when
$d_b (e) = -r \leq 0$,
we have also the following explicit expansion of
$\rho_{b,\varepsilon}(X^e)$ in terms of the distinguished basis
in~$\TT$:
\begin{equation}
\label{eq:rho-explicit}
\rho_{b,\varepsilon}(X^e) = \sum_{p=0}^r
\binom{r}{p}_{q^{d(b)/2}} X^{e + \varepsilon p b},
\end{equation}
where we use the notation
\begin{equation}
\label{eq:centered-binomial-coeff}
\binom{r}{p}_t =
\frac{(t^r - t^{-r}) \cdots (t^{r-p+1} - t^{-r+p-1})}
{(t^p - t^{-p}) \cdots (t - t^{-1})} \ .
\end{equation}
This expansion follows
from the first case in \eqref{eq:rho-extended}
with the help of the well-known ``$t$-binomial formula"
\begin{equation}
\label{eq:centered-binomial}
\prod_{p=0}^{r-1} (1 + t^{r-1-2p} x) =
\sum_{p=0}^r
\binom{r}{p}_t x^p \ .
\end{equation}

\subsection{Toric frames}

\begin{definition}
\label{def:toric-frame}
A \emph{toric frame} in $\FF$ is a mapping $M : \ZZ^m \to \FF - \{0\}$
of the form
\begin{equation}
\label{eq:auto-basis-presentation}
M(c) = \varphi (X^{\eta(c)}),
\end{equation}
where $\varphi$ is an automorphism of $\FF$, and
$\eta: \ZZ^m \to L$ is an isomorphism of lattices.
\end{definition}

Note that both $\varphi$ and $\eta$ are not uniquely determined by
a toric frame~$M$.

By the definition, the elements $M(c)$ form a $\ZZ[q^{\pm 1/2}]$-
basis of an isomorphic copy $\varphi (\TT)$ of the based quantum
torus $\TT$; their multiplication and commutation relations are given by
\begin{equation}
\label{eq:mult-Mc}
M(c) M(d) =  q^{\Lambda_M(c,d)/2} M(c+d),
\end{equation}
and
\begin{equation}
\label{eq:comm-Mc}
M(c) M(d) =  q^{\Lambda_M(c,d)} M(d)M(c),
\end{equation}
where the bilinear form $\Lambda_M$ on $\ZZ^m$ is obtained by
transferring the form $\Lambda$ from $L$ by means of the
lattice isomorphism $\eta$.
(Note that either of \eqref{eq:mult-Mc} and \eqref{eq:comm-Mc} establishes in particular that
$\Lambda_M$ is well defined, i.e., does not depend on the choice of~$\eta$.)
In view of \eqref{eq:unit}, we have
\begin{equation}
\label{eq:M-unit}
M(0) = 1, \quad M(c)^{-1} = M(-c) \quad (c \in \ZZ^m) \, .
\end{equation}
We denote by the same symbol $\Lambda_M$
the corresponding $m \times m$ integer matrix with entries
\begin{equation}
\label{eq:LambdaMij}
\lambda_{ij} = \Lambda_M (e_i, e_j),
\end{equation}
where $\{e_1, \dots, e_m\}$ is the standard basis of $\ZZ^m$.

Given a toric frame, we set $X_i = M(e_i)$ for $i \in [1,m]$.
In view of \eqref{eq:comm-Mc}, the elements
$X_i$ \emph{quasi-commute}:
\begin{equation}
\label{eq:Xi-q-com}
X_i X_j = q^{\lambda_{ij}} X_j X_i \, .
\end{equation}
In the ``classical limit" $q = 1$, the set
$\tilde {\mathbf{X}} = \{X_1, \dots, X_m\}$
specializes to an (arbitrary)
free generating set $\tilde \xx$ of the ambient field, while the set
$\{M(c): c \in \ZZ^m\}$ turns into the set of all Laurent monomials
in the elements of $\tilde \xx$.

\begin{lemma}
\label{lem:frame-positive-part}
A toric frame $M : \ZZ^m \to \FF - \{0\}$
is uniquely determined by the elements
$X_i = M(e_i)$ for $i \in [1,m]$.
\end{lemma}

\begin{proof}
In view of \eqref{eq:mult-Mc}, \eqref{eq:LambdaMij}, and
\eqref{eq:Xi-q-com}, we get
\begin{equation}
\label{eq:M-Xi}
M(a_1, \dots, a_m) = q^{\frac{1}{2} \sum_{\ell < k} a_k a_\ell
\lambda_{k \ell}} X_1^{a_1} \cdots X_m^{a_m}
\end{equation}
for any $(a_1, \dots, a_m) \in \ZZ^m$, which implies our
statement.
\end{proof}

In spite of Lemma~\ref{lem:frame-positive-part}, we still prefer
to include the whole infinite family of elements~$M(c)$ into
Definition~\ref{def:toric-frame}, since there seems to be no nice
way to state the needed conditions in terms of the finite set
$\tilde {\mathbf{X}}$.

\subsection{Quantum seeds and their mutations}

Now everything is ready for a quantum analog of Definition~\ref{def:seed}.

\begin{definition}
\label{def:quantum-seed}
A \emph{quantum seed} is a pair $(M, \tilde B)$,
where
\begin{itemize}
\item $M$ is a toric frame in $\FF$.

\item $\tilde B$ is an $m \times n$ integer matrix
with rows labeled by $[1,m]$ and columns labeled by an
$n$-element subset $\ex \subset [1,m]$.

\item The pair
$(\Lambda_M, \tilde B)$
is compatible in the sense of Definition~\ref{def:compatible-triple}.
\end{itemize}
As in Definition~\ref{def:seed}, quantum seeds are defined up to a permutation of
the standard basis in $\ZZ^m$ together with the corresponding relabeling of rows and columns of $\tilde B$.
\end{definition}

\begin{remark}
In the ``classical limit" $q =1$, the quasi-commutation relations
\eqref{eq:comm-Mc} give rise to the Poisson structure on the
cluster algebra introduced and studied in \cite{gsv}.
In fact, the compatibility condition for the pair $(\Lambda_M, \tilde B)$
appears in \cite[(1.7)]{gsv}.
Furthermore, for $k \in \ex$, let $b^k \in \ZZ^m$ denote the $k$th column
of~$\tilde B$.
As a special case of \eqref{eq:comm-Mc}, for every $j, k \in \ex$, we get
$$M(b^j) M(b^k) =  q^{\Lambda_M(b^j,b^k)} M(b^k)M(b^j),$$
where the exponent $\Lambda_M(b^j,b^k)$ is the $(j,k)$-entry of
the matrix $\tilde B^T \Lambda_M \tilde B$.
Since the pair $(\Lambda_M, \tilde B)$ is compatible,
this exponent is equal to $d_j b_{jk} = - d_k b_{kj}$, where
the positive integers $d_j$ for $j \in \ex$ have the same
meaning as in Definition~\ref{def:compatible-triple}.
In the limit $q = 1$, this agrees with the calculation of the Poisson structure
from \cite[Theorem~1.4]{gsv} in the so-called $\tau$-coordinates.
\end{remark}

Our next target is a quantum analogue of
Definition~\ref{def:seed-mutation}.
Let $(M, \tilde B)$ be a quantum seed.
Fix an index $k \in \ex$ and a sign $\varepsilon \in \{\pm 1\}$.
We define a mapping $M' : \ZZ^m \to \FF - \{0\}$
by setting, for $c = (c_1, \dots, c_m) \in \ZZ^m$ with $c_k \geq 0$,
\begin{equation}
\label{eq:M'}
M'(c) = \sum_{p=0}^{c_k} \binom{c_k}{p}_{q^{d_k/2}}
M(E_\varepsilon c + \varepsilon p b^k),
\quad M'(-c) = M'(c)^{-1},
\end{equation}
where we use the $t$-binomial coefficients from \eqref{eq:centered-binomial-coeff},
the matrix $E_\varepsilon$ is given by \eqref{eq:E-entries},
and the vector $b^k \in \ZZ^m$ is the $k$th column of $\tilde B$.
Finally, let $\tilde B' = \mu_k (\tilde B)$ be given by \eqref{eq:matrix-mutation}.

\begin{proposition}
\label{pr:quantum-seed-mutation}
\begin{enumerate}
\item The mapping $M'$ is a toric frame independent of the choice
of a sign~$\varepsilon$.

\item The pair $(\Lambda_{M'}, \tilde B')$ is obtained from
$(\Lambda_M, \tilde B)$ by the mutation in direction~$k$
(see Definition~\ref{def:triple-mutation}).

\item The pair $(M', \tilde B')$ is a quantum seed.
\end{enumerate}
\end{proposition}

\proof (1) To  see that $M'$ is independent of the choice
of~$\varepsilon$, notice that the summation term in
\eqref{eq:M'} does not change if we replace~$\varepsilon$
with~$-\varepsilon$, and $p$ with $c_k - p$ (this is a
straightforward check).
To show that $M'$ is a toric frame, we express $M$ according to
\eqref{eq:auto-basis-presentation}.
Replacing the initial based quantum torus $\TT$ with
$\varphi (\TT)$, and using $\eta$ to identify the lattice
$L$ with $\ZZ^m$, we may assume from the start that $L = \ZZ^m$,
and $M(c) = X^c$ for any $c \in L$.
Note that the compatibility condition for the pair
$(\Lambda_M, \tilde B)$ can be simply written as
\begin{equation}
\label{eq:bk-properties}
\Lambda (b^j, e_i) = \delta_{ij} d_j \quad (i \in [1,m], \,
j \in \ex) \, .
\end{equation}
It follows that, using the notation introduced in
Section~\ref{sec:auto}, we get $d(b^k) = d_k$ for
$k \in \ex$, and
$d_{b^k} (E_\varepsilon c) = - c_k$.
Comparing \eqref{eq:M'} with \eqref{eq:rho-explicit}, we now obtain
\begin{equation}
\label{M'-right-way}
M'(c) = \rho_{b^k,\varepsilon} (X^{E_\varepsilon c})  \quad (c \in L);
\end{equation}
thus, $M'$ is of the form \eqref{eq:auto-basis-presentation},
i.e., is a toric frame.

\medskip

(2) In view of \eqref{eq:LambdaMij} and \eqref{M'-right-way},
the matrices $\Lambda_{M'}$ and $\Lambda_{M}$ are related by
$\Lambda_{M'} = E_\varepsilon^T \Lambda_{M} E_\varepsilon$,
so the claim follows from \eqref{eq:Lambda-mutated}.

\medskip

(3) The statement follows from parts (1) and (2) in view of
Proposition~\ref{pr:triple-mutation}.
\endproof

Proposition~\ref{pr:quantum-seed-mutation} justifies the following definition.

\begin{definition}
\label{def:quantum-seed-mutation}
Let $(M, \tilde B)$ be a quantum seed, and $k \in \ex$.
We say that the quantum seed $(M', \tilde B')$ given by \eqref{eq:M'}
and \eqref{eq:matrix-mutation} is obtained from $(M, \tilde B)$
by the \emph{mutation} in direction~$k$, and write
$(M', \tilde B') = \mu_k (M, \tilde B)$.
\end{definition}

The following proposition demonstrates that Definition~\ref{def:quantum-seed-mutation}
is indeed a quantum analogue of Definition~\ref{def:seed-mutation}.

\begin{proposition}
\label{pr:X-to-X'}
Let $(M, \tilde B)$ be a quantum seed, and
suppose the quantum seed $(M', \tilde B')$
is obtained from $(M, \tilde B)$ by the mutation
in direction~$k \in \ex$.
For $i \in [1,m]$, let $X_i = M(e_i)$ and $X'_i = M'(e_i)$.
Then $X'_i = X_i$ for $i \neq k$, and
$X'_k$ is given by the following quantum analogue of
the exchange relation {\rm \eqref{eq:exchange-rel-xx}}:
\begin{equation}
\label{eq:quantum-exchange}
X'_k = M(- e_k + \sum_{b_{ik} > 0} b_{ik} e_i) \ + \
 M(- e_k - \sum_{b_{ik} < 0} b_{ik} e_i) \, .
\end{equation}
\end{proposition}

\proof This follows at once by applying
\eqref{eq:M'} to $c = e_i$ for $i \in [1,m]$.
\endproof


\begin{proposition}
\label{pr:quantum-mutation-involutive}
The mutation of quantum seeds is involutive: if
$(M', \tilde B') = \mu_k (M, \tilde B)$
then $\mu_k (M', \tilde B') = (M, \tilde B)$.
\end{proposition}

\begin{proof}
As in the proof of Proposition~\ref{pr:quantum-seed-mutation},
we can assume without loss of generality that $L = \ZZ^m$,
and $M(c) = X^c$ for any $c \in L$.
Then the toric frame $M'$ is given by \eqref{M'-right-way}.
Applying \eqref{M'-right-way} once again, with $\varepsilon$ replaced by
$- \varepsilon$, we see that the toric frame $M''$ in the quantum seed
$\mu_k (M', \tilde B')$ is given by
$$M''(c) =  \rho_{b^k,\varepsilon} \rho_{- E_\varepsilon b^k, - \varepsilon}
(X^{E_\varepsilon E'_{-\varepsilon} c}) \ ,$$
where the matrix $E'_{-\varepsilon}$ is given by \eqref{eq:E-entries} applied to
$\tilde B'$ instead of $\tilde B$.
Using an obvious fact that $E_\varepsilon b^k = b^k$
together with \eqref{E'-E}, \eqref{eq:EF-squares}, and
\eqref{eq:rho-mult}, we conclude that $M''(c) = X^c = M(c)$, as required.
\end{proof}

\subsection{Quantum cluster algebras}
\label{sec:quantum-ca-def}
In view of Proposition~\ref{pr:quantum-mutation-involutive},
the following relation on quantum seeds is an equivalence
relation: we say that two quantum seeds are
\emph{mutation-equivalent} if they can be obtained from each other
by a sequence of quantum seed mutations.
For a quantum seed $(M, \tilde B)$, we denote by
$\tilde {\mathbf{X}}= \{X_1, \dots, X_m\}$ the corresponding
``free generating set" in $\FF$ given by $X_i = M(e_i)$.
As for the ordinary seeds, we call the subset
${\mathbf{X}}= \{X_j : j \in \ex\} \subset \tilde {\mathbf{X}}$
the \emph{cluster} of the quantum seed $(M, \tilde B)$, and set
$\mathbf{C} = \tilde {\mathbf{X}} - \mathbf{X}$.
The following result is an immediate consequence of
Proposition~\ref{pr:X-to-X'}.

\begin{proposition}
\label{pr:C}
The $(m-n)$-element set $\mathbf{C} = \tilde {\mathbf{X}} - \mathbf{X}$
depends only on the mutation-equivalence class of
a quantum seed $(M, \tilde B)$.
\end{proposition}

Now everything is in place for defining quantum cluster algebras.

\begin{definition}
\label{def:quantum-cluster-algebra}
Let $\mathcal{S}$ be a mutation-equivalence class of quantum seeds in
$\FF$, and let $\mathbf{C} \subset \FF$ be the $(m-n)$-element set
associated to $\mathcal{S}$ as in Proposition~\ref{pr:C}.
The cluster algebra $\AA(\mathcal{S})$ associated with $\mathcal{S}$
is the $\ZZ[q^{\pm 1/2}]$-subalgebra of the ambient skew-field $\FF$,
generated by the union of clusters of all seeds in $\mathcal{S}$,
together with the elements of $\mathbf{C}$ and their inverses.
\end{definition}

Since $\mathcal{S}$ is uniquely determined by each of its quantum seeds
$(M, \tilde B)$, we sometimes denote $\AA(\mathcal{S})$ as $\AA(M,\tilde B)$,
or even simply $\AA(\Lambda_M, \tilde B)$, because a compatible matrix pair
$(\Lambda_M, \tilde B)$ determines this algebra uniquely up to
an automorphism of the ambient skew-field~$\FF$.
We denote by $\PP$ the multiplicative group generated by $q^{1/2}$
and $\mathbf{C}$, and treat the integer group ring $\ZZ \PP$ as the
\emph{ground ring} for the cluster algebra.
In other words, $\ZZ \PP$ is the ring of
Laurent polynomials in the elements of $\mathbf{C}$ with
coefficients in $\ZZ[q^{\pm 1/2}]$.

\section{Upper bounds and quantum Laurent phenomenon}
\label{sec:upper-bound}

Let $(M, \tilde B)$ be a quantum seed in~$\FF$, and
$\tilde {\mathbf{X}}= \{X_1, \dots, X_m\}$ denote the corresponding
``free generating set" in $\FF$ given by $X_i = M(e_i)$.
As in \cite{bfz-clust3}, we will associate with $(M, \tilde B)$ a subalgebra
$\mathcal{U}(M, \tilde B) \subset \FF$ called the (quantum)
\emph{upper cluster algebra}, or simply the \emph{upper bound}.

Let $\ZZ \PP[\mathbf{X}^{\pm 1}]$ denote the based
quantum torus generated by $\tilde {\mathbf{X}}$; this is a
$\ZZ[q^{\pm 1/2}]$-subalgebra of $\FF$ with the basis
$\{M(c) : c \in \ZZ^m\}$.
For the sake of convenience, in this section we assume that
$\tilde {\mathbf{X}}$ is numbered so that
its cluster $\mathbf{X}$ has the form
$\mathbf{X} = \{X_1, \dots, X_n\}$.
Thus, the complement
$\mathbf{C} =  \tilde {\mathbf{X}} - \mathbf{X}$ is given by
$\mathbf{C} = \{X_{n+1}, \dots, X_m\}$, and
the ground ring $\ZZ \PP$ is the ring of integer Laurent
polynomials in the (quasi-commuting) variables $q^{1/2}, X_{n+1}, \dots, X_m$.
For $k \in [1,n]$, let $(M_k, \tilde B_k)$
denote the quantum seed obtained from $(M, \tilde B)$
by the mutation in direction~$k$, and let
${\mathbf{X}}_k$ denote its cluster; thus, we have
\begin{equation}
\label{eq:quantum-exchange-clusters}
\mathbf{X}_k = \mathbf{X} - \{X_k\} \cup \{X'_k\},
\end{equation}
where $X'_k$ is given by \eqref{eq:quantum-exchange}.

Following \cite[Definition~1.1]{bfz-clust3}, we denote by
$\mathcal{U}(M, \tilde B) \subset \FF$
the $\ZZ \PP$-subalgebra of $\FF$ given by
\begin{equation}
\label{eq:upper-bound}
\mathcal{U}(M, \tilde B) =
\ZZ \PP[\mathbf{X}^{\pm 1}]
\cap \ZZ \PP[\mathbf{X}_1^{\pm 1}]
\cap \cdots \cap \ZZ \PP[\mathbf{X}_n^{\pm 1}]\, .
\end{equation}
In other words, $\mathcal{U}(M, \tilde B)$ is formed by the elements of $\FF$
which are expressed as Laurent polynomials over $\ZZ \PP$ in the variables from
each of the clusters $\mathbf{X}, \mathbf{X}_1, \dots, \mathbf{X}_n$.

Our first main result is a quantum analog of
\cite[Theorem~1.5]{bfz-clust3}.

\begin{theorem}
\label{th:upper-bound-invariance}
The algebra $\mathcal{U}(M, \tilde B)$
depends only on the mutation-equivalence class of the quantum seed $(M, \tilde B)$.
\end{theorem}

Theorem~\ref{th:upper-bound-invariance} justifies the notation
$\mathcal{U}(M, \tilde B) =
\mathcal{U}(\mathcal{S})$, where $\mathcal{S}$ is
the mutation-equivalence class of $(M, \tilde B)$;
in fact, we have
\begin{equation}
\label{eq:US-intersection}
\mathcal{U}(\mathcal{S}) =
\bigcap_{(M, \tilde B) \in \mathcal{S}}
\ZZ \PP[\mathbf{X}^{\pm 1}] \ .
\end{equation}
In view of Propositions~\ref{pr:X-to-X'} and \ref{pr:quantum-mutation-involutive},
$\tilde {\mathbf{X}} \subset \mathcal{U}(\mathcal{S})$ for every
quantum seed $(M, \tilde B)$ in~$\mathcal{S}$.
Therefore, Theorem~\ref{th:upper-bound-invariance} has the
following important corollary that justifies calling $\mathcal{U}(\mathcal{S})$
the \emph{upper bound} for the cluster algebra.

\begin{corollary}
\label{cor:A-subset-U}
The cluster algebra $\AA(\mathcal{S})$ is contained in
$\mathcal{U}(\mathcal{S})$.
Equivalently, $\AA(\mathcal{S})$ is contained in
the quantum torus
$\ZZ \PP[\mathbf{X}^{\pm 1}]$
for every quantum seed
$(M, \tilde B) \in \mathcal{S}$ with the cluster $\mathbf{X}$
(we refer to this property as the quantum Laurent phenomenon).
\end{corollary}

\begin{example}
\label{ex:quantum-rank-2}
Let $\AA(b,c)$ be the quantum cluster algebra associated with a
compatible pair $(\Lambda, \tilde B)$ of the form
$$\Lambda = \mat{0}{1}{-1}{0}, \quad
\tilde B = B = \mat{0}{b}{-c}{0}$$
for some positive integers~$b$ and~$c$.
Tracing the definitions, we see that $\AA(b,c)$ can be described
as follows (cf.~\cite{fz-clust1,sz}).
The ambient field $\FF$ is the skew-field of fractions of the quantum torus
with generators $Y_1$ and $Y_2$ satisfying the quasi-commutation relation
$Y_1 Y_2 = q Y_2 Y_1$.
Then $\AA(b,c)$ is the $\ZZ[q^{\pm 1/2}]$-subalgebra of~$\FF$
generated by a sequence of cluster variables $\{Y_m: m \in \ZZ\}$ defined
recursively from the relations
\begin{equation}
\label{eq:clusterrelations}
Y_{m-1} Y_{m+1} =
\left\{
\begin{array}[h]{ll}
q^{b/2} Y_m^b + 1 & \quad \mbox{$m$ odd;} \\
q^{c/2} Y_m^c + 1 & \quad \mbox{$m$ even.}
\end{array}
\right.
\end{equation}
The clusters are the pairs $\{Y_m, Y_{m+1}\}$ for all $m \in \ZZ$.
One checks easily that
$$Y_m Y_{m+1} = q Y_{m+1} Y_m \quad (m \in \ZZ).$$
According to Corollary~\ref{cor:A-subset-U}, every cluster
variable~$Y_m$ is a Laurent polynomial in $Y_1$ and $Y_2$ with
coefficients in $\ZZ[q^{\pm 1/2}]$.
A direct calculation gives these polynomials explicitly in the
\emph{finite type} cases when $bc \leq 3$
(cf.~\cite[(4.4)-(4.6)]{sz}).
In accordance with \eqref{eq:M-Xi}, in the following formulas we
use the notation
$$Y^{(a_1, a_2)} = q^{-a_1 a_2/2} Y_1^{a_1} Y_2^{a_2} \quad
(a_1, a_2 \in \ZZ).$$

\smallskip

\noindent  \textsl{Type $A_2$:} $(b,c) = (1,1)$.
\begin{eqnarray}
\label{eq:Ym-Laurent-A2}
&&Y_3 = Y^{(-1,1)} + Y^{(-1,0)}, \quad Y_4 =
 Y^{(0,-1)} + Y^{(-1,-1)} + Y^{(-1,0)},\\
\nonumber
&&Y_5 = Y^{(1,-1)} + Y^{(0,-1)}, \quad Y_6 = Y_1, \quad Y_7 = Y_2 \ .
\end{eqnarray}

\smallskip

\noindent \textsl{Type $B_2$:} $(b,c) = (1,2)$.
\begin{eqnarray}
\label{eq:Ym-Laurent-B2}
&&Y_3 = Y^{(-1,2)} + Y^{(-1,0)}, \quad
Y_4 = Y^{(0,-1)} + Y^{(-1,-1)} + Y^{(-1,1)},\\
\nonumber
&&Y_5 = Y^{(1,-2)} + (q^{1/2} + q^{-1/2}) Y^{(0,-2)}
+ Y^{(-1,-2)} + Y^{(-1,0)},\\
\nonumber
&&Y_6 = Y^{(1,-1)} + Y^{(0,-1)}, \quad Y_7 = Y_1, \quad Y_8 = Y_2 \ .
\end{eqnarray}

\smallskip

\noindent  \textsl{Type $G_2$:} $(b,c) = (1,3)$.
\begin{eqnarray}
\label{eq:Ym-Laurent-G2}
&Y_3 =& Y^{(-1,3)} + Y^{(-1,0)}, \quad
Y_4 = Y^{(0,-1)} + Y^{(-1,-1)} + Y^{(-1,2)},\\
\nonumber
&Y_5=& Y^{(1,-3)} + (q+1+q^{-1})(Y^{(0,-3)} + Y^{(-1,0)}
+ Y^{(-1,-3)})\\
\nonumber
&&+  Y^{(-2,3)} + (q^{3/2} + q^{-3/2}) Y^{(-2,0)} + Y^{(-2,-3)},\\
\nonumber
&Y_6=& Y^{(1,-2)} + (q^{1/2} + q^{-1/2}) Y^{(0,-2)} + Y^{(-1,-2)}
+ Y^{(-1,1)},\\
\nonumber
&Y_7 =& Y^{(2,-3)} + (q+1+q^{-1})(Y^{(1,-3)} + Y^{(0,-3)})
+ Y^{(-1,-3)} + Y^{(-1,0)},\\
\nonumber
&Y_8 =& Y^{(1,-1)} + Y^{(0,-1)}, \quad Y_9 = Y_1, \quad Y_{10} = Y_2 \ .
\end{eqnarray}
\end{example}

\smallskip

The rest of this section is devoted to the proof of
Theorem~\ref{th:upper-bound-invariance}.
The proof follows that of \cite[Theorem~1.5]{bfz-clust3}
but we have to deal with some technical complications caused by non-commutativity
of a quantum torus.
As a rule, the arguments in \cite{bfz-clust3} will require only
obvious changes if the quantum analogs of all participating elements
quasi-commute with each other.
We shall provide more details when more serious changes will be
needed.

We start with an analog of \cite[Lemma~4.1]{bfz-clust3}.

\begin{lemma}
\label{lem:upper-bound-intersection-1}
The algebra $\mathcal{U}(M, \tilde B)$
can be expressed as follows:
\begin{equation}
\label{eq:upper-bound-2}
\mathcal{U}(M, \tilde B) =
\bigcap_{k = 1}^n
\ZZ \PP[X_1^{\pm 1}, \ldots, X_{k-1}^{\pm 1}, X_k, X'_k,
 X_{k+1}^{\pm 1}, \ldots, X_{n}^{\pm 1}],
\end{equation}
where $X'_k$ is given by {\rm \eqref{eq:quantum-exchange}.}
\end{lemma}

\proof
In view of \eqref{eq:upper-bound}, it is enough to show that
\begin{equation}
\label{eq:two-laurent-conditions-2}
\ZZ \PP[\mathbf{X}^{\pm 1}]
\cap \ZZ \PP[\mathbf{X}_1^{\pm 1}] =
\ZZ \PP[X_1,X'_1, X_2^{\pm 1}, \dots, X_n^{\pm 1}].
\end{equation}

As in \cite{bfz-clust3}, \eqref{eq:two-laurent-conditions-2} is a
consequence of the following easily verified properties.

\begin{lemma}
\label{lem:two-laurent}
\begin{enumerate}
\item
Every element $Y \in \ZZ \PP[\mathbf{X}^{\pm 1}]$ can be
uniquely written in the form
\begin{equation}
\label{eq:usual-Laurent}
Y = \sum_{r \in \ZZ} c_r X_1^r\ ,
\end{equation}
where each coefficient $c_r$ belongs to
$\ZZ \PP[X_2^{\pm 1}, \dots, X_n^{\pm 1}]$, and all but
finitely many of them are equal to~$0$.

\item
Every element
$Y \in \ZZ \PP[\mathbf{X}^{\pm 1}]
\cap \ZZ \PP[\mathbf{X}_1^{\pm 1}]$ can be
uniquely written in the form
\begin{equation}
\label{eq:usual-two-Laurent}
Y = c_0 + \sum_{r \geq 1}
(c_r X_1^r + c'_r (X'_1)^r)\ ,
\end{equation}
where all coefficients $c_r$ and $c'_r$ belong to
$\ZZ \PP[X_2^{\pm 1}, \dots, X_n^{\pm 1}]$, and all but
finitely many of them are equal to~$0$.
\end{enumerate}
\end{lemma}

Our next target is an analog of \cite[Lemma~4.2]{bfz-clust3}.
As in the proof of Proposition~\ref{pr:quantum-seed-mutation},
in what follows, we will assume without loss of generality that
$L = \ZZ^m$, and the toric frame of the initial quantum seed $(M, \tilde B)$
is given by $M(c) = X^c$ for any $c \in L$.
In particular, we view the columns $b^j$ of $\tilde B$ as elements of~$L$.
According to \eqref{eq:P-coeff}, for every nonnegative integer $r$ and every
sign $\varepsilon$, we have a well-defined element
$P^r_{b^1,\varepsilon} \in \ZZ \PP[X_2^{\pm 1}, \dots, X_m^{\pm 1}]$.
Note that, in view of \eqref{eq:multiplication-F-com} and
\eqref{eq:bk-properties}, $P^r_{b^1,\varepsilon}$
belongs to the center of the algebra
$\ZZ \PP[X_2^{\pm 1}, \dots, X_m^{\pm 1}]$.
In particular, $P^r_{b^1,+}$ and $P^r_{b^1,-}$
commute with each other; an easy check shows that their ratio
is an invertible element of the center of
$\ZZ \PP[X_2^{\pm 1}, \dots, X_m^{\pm 1}]$.

\begin{lemma}
\label{lem:two-laurent-conditions-1}
An element $Y \in \FF$ belongs to
$\ZZ \PP[X_1,X'_1, X_2^{\pm 1}, \dots, X_n^{\pm 1}]$
if and only if it has the form {\rm \eqref{eq:usual-Laurent}},
and for each $r > 0$, the coefficient $c_{-r}$ is divisible by
$P^r_{b^1,+}$ in the algebra
$\ZZ \PP[X_2^{\pm 1}, \dots, X_n^{\pm 1}]$.
\end{lemma}

\proof  In view of \eqref{M'-right-way} and
\eqref{eq:rho-extended}, we have
\begin{equation}
\label{eq:X'1-powers}
(X'_1)^r = P^r_{b^1,+}
(X^{e'_1})^r \ ,
\end{equation}
where
\begin{equation}
\label{eq:e'1}
e'_1 = - e_1 - \sum_{b_{i1} < 0} b_{i1} e_i \ .
\end{equation}
Combining \eqref{eq:X'1-powers} with \eqref{eq:usual-two-Laurent},
we obtain the desired claim.
\endproof

Our next step is an analog of \cite[Proposition~4.3]{bfz-clust3}.

\begin{proposition}
\label{pr:upper-bound-intersection-2}
Suppose that $n \geq 2$.
Then
\begin{equation}
\label{eq:upper-bound-intersection-2}
\mathcal{U} (M, \tilde B) =
\bigcap_{j = 2}^n
\ZZ \PP[X_1, X'_1, X_2^{\pm 1}, \ldots, X_{j-1}^{\pm 1}, X_j, X'_j,
 X_{j+1}^{\pm 1}, \ldots, X_{n}^{\pm 1}].
\end{equation}
\end{proposition}

\begin{proof}
As in the proof of \cite[Proposition~4.3]{bfz-clust3}, we can
assume that $n = 2$, i.e., the ground ring $\ZZ \PP$ is the ring of Laurent
polynomials in $q, X_3, \dots, X_m$.
Thus, it suffices to show the following
analog of \cite[(4.4)]{bfz-clust3}:
\begin{equation}
\label{eq:rk-2-upper=lower}
\ZZ \PP[X_1,X'_1, X_2^{\pm 1}] \cap
\ZZ \PP[X_1^{\pm 1}, X_2,X'_2] =
\ZZ \PP[X_1,X'_1, X_2,X'_2].
\end{equation}

The proof of (\ref{eq:rk-2-upper=lower}) breaks into two cases.

\medskip

\noindent \emph{Case 1:} $b_{12} = b_{21} = 0$.
In this case, the elements $P^r_{b^1,+}$ and
$P^s_{b^2,+}$ belong to the center of
$\ZZ \PP$ for all $r, s > 0$;
furthermore, $P^r_{b^1,+}$ commutes with $X_2$, while
$P^s_{b^2,+}$ commutes with $X_1$.
Arguing as in \cite{bfz-clust3}, we reduce the proof to the
following statement: if an element of
$\ZZ \PP$ is divisible by each of the $P^r_{b^1,+}$ and
$P^s_{b^2,+}$ then it is divisible by their product.
By Proposition~\ref{pr:center} below, it suffices to check that
$P^r_{b^1,+}$ and $P^s_{b^2,+}$ are relatively
prime in the center of $\ZZ \PP$.
This follows from the fact that $\tilde B$ has full rank (see
Proposition~\ref{pr:stronger}), and so the columns
$b^1$ and $b^2$ are not proportional to each other.

\medskip

\noindent \emph{Case 2:} $b_{12}b_{21} < 0$.
In this case, the proof goes through the same steps as in
\cite{bfz-clust3}, with some obvious modifications taking into
account non-commutativity.
We leave the details to the reader.
\end{proof}

To finish the proof of Theorem~\ref{th:upper-bound-invariance},
it is enough to show that
$\mathcal{U}(M, \tilde B)$ does not
change under the mutation in direction~$1$.
If $n=1$, there is nothing to prove, so we assume that $n \geq 2$.
Let $X''_2$ be the cluster variable that replaces $X_2$ in the
cluster $\mathbf{X}_1$ under the mutation in direction~$2$.
In view of \eqref{eq:upper-bound-intersection-2},
Theorem~\ref{th:upper-bound-invariance} becomes a consequence
of the following lemma.

\begin{lemma}
\label{lem:reduction-to-rank-2}
In the above notation, we have
$$\ZZ \PP [X_1,X'_1, X_2,X'_2]
= \ZZ \PP[X_1,X'_1, X_2,X''_2].$$
\end{lemma}

\begin{proof}
By symmetry, it is enough to show that
\begin{equation}
\label{eq:X2''}
X''_2 \in \ZZ \PP[X_1,X'_1, X_2,X'_2]\ .
\end{equation}
The following proof of \eqref{eq:X2''} uses the same strategy as in the
proof of \cite[Lemma~4.6]{bfz-clust3}, but one has to keep a
careful eye on the non-commutativity effects.

We start by recalling the assumption that $L = \ZZ^m$, and the
initial toric frame $M$ is given by $M(c) = X^c$ for any $c \in L$.
Then the toric frames of the adjacent quantum seeds
are given by \eqref{M'-right-way}.
For typographic reasons, we rename the quantum seed
$(M_1, \tilde B_1)= \mu_1(M, \tilde B)$ to $(M', \tilde B')$
(so the entries of the matrix $\tilde B_1 = \tilde B'$
are denoted $b'_{ij}$), and also use the notation
$(M'', \tilde B'')= \mu_2(M', \tilde B')$.
Thus, $X''_2 = M''(e_2)$.
Without loss of generality, we assume that the matrix entry
$b_{12}$ of $\tilde B$ is non-positive; and we set
$r = - b_{12} \geq 0$.
Since the principal parts of $\tilde B$ and $\tilde B'$ are
skew-symmetrizable,
it follows that $b_{21} \geq 0$, $b'_{12} = r$, and $b'_{21} \leq 0$.

Applying \eqref{eq:quantum-exchange} and \eqref{M'-right-way},
we see that
$$X''_2 = M'(e''_2) + M'(e''_2 + (b')^2)
= \rho_{b^1,+} (X^{E_+ e''_2} + X^{E_+ (e''_2 + (b')^2)})
\ ,$$
where
\begin{equation}
\label{eq:e''2}
e''_2 = - e_2 - \sum_{i>2, \ b'_{i2} < 0} b'_{i2} e_i \, ,
\end{equation}
$(b')^2$ is the second column of $\tilde B'$, and
$E_+$ is given by \eqref{eq:E-entries} with $k = 1$.
Note that the summation in \eqref{eq:e''2} does not include a multiple of
$e_1$ because $b'_{12} = r \geq 0$; this implies that
$E_+ e''_2 = e''_2$.
We also have $E_+ (b')^2 = b^2$ (to see this, use \eqref{eq:mutation-product}
to write $\tilde B' = E_+ \,\tilde B\, F_+$,
and note that the second column of $\tilde B \, F_+$ is equal to
$b^2$, hence $(b')^2 = E_+ b^2$, and so our statement follows from
\eqref{eq:EF-squares}).
Remembering \eqref{eq:rho-extended} and \eqref{eq:bk-properties}, we conclude that
\begin{equation}
\label{eq:X''2-expression}
X''_2 = X^{e''_2} + P_{b^1,+}^r X^{e''_2 + b^2} \ .
\end{equation}

On the other hand, setting
$$e'_2 = - e_2 - \sum_{b_{i2} < 0} b_{i2} e_i \ ,$$
we have
$$X'_2 = X^{e'_2} + X^{e'_2 + b^2} \, ;$$
applying \eqref{eq:multiplication-F} and \eqref{eq:bk-properties},
we obtain
\begin{equation}
\label{eq:X2X'2}
q^{- \Lambda (e_2, e'_2)/2} X_2 X'_2
= X^{e_2 + e'_2} + q^{- d_2/2} X^{e_2 + e'_2 + b^2} \ .
\end{equation}
Note that the second summand $F = q^{- d_2/2} X^{e_2 + e'_2 + b^2}$
is an invertible element of $\ZZ \PP$;
thus, to prove the desired inclusion \eqref{eq:X2''},
it suffices to show that
$$X''_2 F \in \ZZ \PP[X_1,X'_1, X_2,X'_2]\ .$$

Using \eqref{eq:X''2-expression} and \eqref{eq:X2X'2}, we write
$$X''_2 F = q^{- \Lambda (e_2, e'_2)/2} S_1 - S_2 + S_3,$$
where
\begin{align*}
&S_1 =  P_{b^1,+}^r X^{e''_2 + b^2}  X_2 X'_2, \\[.1in]
& S_2 = (P_{b^1,+}^r - 1)  X^{e''_2 + b^2} X^{e_2 + e'_2}, \\[.1in]
& S_3 = q^{- d_2/2} X^{e''_2} X^{e_2 + e'_2 + b^2} -
X^{e''_2 + b^2} X^{e_2 + e'_2}.
\end{align*}
To complete the proof, we will show that
$$S_1, S_2 \in \ZZ \PP[X_1,X'_1, X_2,X'_2],
\quad S_3 = 0.$$

First, we use \eqref{eq:X'1-powers} to rewrite $S_1$ as
\begin{equation}
\label{eq:S1}
S_1 = (X'_1)^r (X^{e'_1})^{-r} X^{e''_2 + b^2} X_2 X'_2 \ .
\end{equation}
A direct check shows that the vector $-r e'_1 + e''_2 + b^2 + e_2$
has the first two components equal to~$0$; it follows that the middle
factor $(X^{e'_1})^{-r} X^{e''_2 + b^2} X_2$ in \eqref{eq:S1}
is an invertible element of $\ZZ \PP$.
Thus, $S_1 \in \ZZ \PP[X_1,X'_1, X_2,X'_2]$, as desired.

To show the same inclusion for $S_2$, we notice
that $P_{b^1,+}^r - 1$ is a polynomial in $X^{b^1}$ with
coefficients in $\ZZ[q^{\pm 1/2}]$ and zero constant term.
If $r = - b_{12} = 0$ then $S_2 = 0$, and there is nothing to prove.
Otherwise, the desired inclusion follows from the fact
that the first two components of $b^1$ are $(0, b_{21})$ with
$b_{21} > 0$, while the first two components of
$e''_2 + b^2 + e_2 + e'_2$ are $(0,-1)$.

Finally, to show that $S_3 = 0$, in view of \eqref{eq:multiplication-F},
we only need to check that
$$- d_2 + \Lambda(e''_2, e_2 + e'_2 + b^2) =
\Lambda (e''_2 + b^2, e_2 + e'_2),$$
or, equivalently,
$$\Lambda (b^2, e_2 + e'_2 + e''_2) = -d_2,$$
which is a direct consequence of \eqref{eq:bk-properties}.
This completes the proof of Lemma~\ref{lem:reduction-to-rank-2} and
Theorem~\ref{th:upper-bound-invariance}.
\end{proof}

\section{Exchange graphs, bar-involutions, and gradings}
\label{sec:exchange-graphs}

Recall that the \emph{exchange graph} of the cluster algebra
$\AA(\mathcal{S})$ associated with a mutation-equivalent
class of seeds $\mathcal{S}$ has the seeds from $\mathcal{S}$ as vertices,
and the edges corresponding to seed mutations
(cf.~\cite[Section~7]{fz-clust1} or \cite[Section~1.2]{fz-clust2}).
We define the exchange graph of a quantum cluster algebra in
exactly the same way: the vertices correspond to its quantum seeds,
and the edges to quantum seed mutations.
As explained in Section~\ref{sec:quantum-ca-def}, we can associate
the quantum cluster algebra with a compatible matrix pair $(\Lambda_M, \tilde B)$,
and denote it $\AA(\Lambda_M, \tilde B)$.
Let $E(\Lambda_M, \tilde B)$ denote the exchange graph of
$\AA(\Lambda_M, \tilde B)$, and $E(\tilde B)$ denote the exchange graph of
the cluster algebra $\AA(\tilde B)$ obtained from $\AA(\Lambda_M, \tilde B)$
by the specialization $q = 1$.
Then the graph $E(\Lambda_M, \tilde B)$ naturally covers $E(\tilde B)$.

\begin{theorem}
\label{th:exchange-graphs}
The specialization $q = 1$ identifies the quantum exchange graph
$E(\Lambda_M, \tilde B)$ with the ``classical"
exchange graph $E(\tilde B)$.
\end{theorem}

The proof of Theorem~\ref{th:exchange-graphs} will require a little preparation.
For a quantum seed $(M,\tilde B)$, let $\TT_M$ denote the
corresponding based quantum torus having $\{M(c): c \in \ZZ^m\}$
as a $\ZZ[q^{\pm 1/2}]$-basis.
This is the same algebra that was previously denoted by
$\ZZ \PP[\mathbf{X}^{\pm 1}]$, where $\mathbf{X}$ is the cluster
of $(M,\tilde B)$; thus, we can rewrite \eqref{eq:US-intersection} as
\begin{equation}
\label{eq:US-intersection-2}
\mathcal{U}(\mathcal{S}) =
\bigcap_{(M, \tilde B) \in \mathcal{S}} \TT_M \ ,
\end{equation}
where $\mathcal{S}$ is the mutation-equivalence class of $(M,\tilde B)$.
We associate with $(M,\tilde B)$ the $\ZZ$-linear \emph{bar-involution}
$X \mapsto \overline X$ on $\TT_M$ by setting
\begin{equation}
\label{eq:bar-involution-trivial}
\overline {q^{r/2} M(c)} = q^{-r/2} M(c) \quad (r \in \ZZ, \, c \in \ZZ^m)\ .
\end{equation}

\begin{proposition}
\label{pr:bar-involution-compatible}
Let $\mathcal{S}$ be the mutation-equivalence class of
a quantum seed $(M,\tilde B)$.
Then the bar-involution associated with $(M,\tilde B)$ preserves
the subalgebra $\mathcal{U}(\mathcal{S}) \subset \TT_M$,
and its restriction to $\mathcal{U}(\mathcal{S})$ depends only
on~$\mathcal{S}$.
\end{proposition}

\begin{proof}
It suffices to show the following: if two quantum seeds
$(M, \tilde B)$ and $(M', \tilde B')$
are obtained from each other by a mutation in some direction~$k$, then
the corresponding bar-involutions have the same restriction to
$\TT_M \cap \TT_{M'}$.
Using \eqref{eq:usual-two-Laurent}, we see that each element
of $\TT_M \cap \TT_{M'}$ is a $\ZZ[q^{\pm 1/2}]$-linear
combination of the elements $M(c)$ and $M'(c)$ for all
$c \in \ZZ^m$ with $c_k \geq 0$.
It remains to observe that, in view of \eqref{eq:M'},
each $M'(c)$ with $c_k \geq 0$ is invariant under the
bar-involution associated with $(M, \tilde B)$.
\end{proof}

\noindent {\bf Proof of Theorem~\ref{th:exchange-graphs}.}
We need to show the following: if two quantum seeds
$(M, \tilde B)$ and $(M', \tilde B')$
are mutation-equivalent, and such that $\tilde B' = \tilde B$
and $M'(c)|_{q=1} = M(c)|_{q=1}$ for all $c \in \ZZ^m$, then~$M' = M$.
(Recall that a quantum seed is defined up to a permutation of the
coordinates in $\ZZ^m$ together with the corresponding relabeling
of rows and columns of $\tilde B$.)
In view of Lemma~\ref{lem:frame-positive-part}, it suffices to
show that $M'(c) = M(c)$ for $c$ being one of the standard basis vectors
$e_1, \dots, e_n$.

By Corollary~\ref{cor:A-subset-U}, $M'(c) \in \TT_M$, i.e.,
$M'(c)$ is a $\ZZ[q^{\pm 1/2}]$-linear combination of the elements
$M(d)$ for $d \in \ZZ^m$.
Let $N(c)$ denote the \emph{Newton polytope} of $M'(c)$, i.e., the
convex hull in $\RR^m$ of the set of all $d \in \ZZ^m$ such that
$M(d)$ occurs in $M'(c)$ with a non-zero coefficient.
We claim that $N(c)$ does not shrink under the specialization $q = 1$,
i.e., that none of the coefficients at vertices of $N(c)$ vanish
under this specialization.
To see this, note that, in view of \eqref{eq:M'},
$M'(c)$ is obtained from a family
$\{M(d) \ : \ d \in \ZZ^m\}$ by a sequence of subtraction-free
rational transformations.
This implies in particular that, whenever~$d$ is a vertex of
$N(c)$, the coefficient of $M(d)$ in $M'(c)$ is a Laurent
polynomial in~$q^{1/2}$ which can also be written as a subtraction-free
rational expression.
Therefore, this coefficient does not vanish at $q=1$, as claimed.
This allows us to conclude that the assumption $M'(c)|_{q=1} = M(c)|_{q=1}$
implies that $M'(c) = p \ M(c)$ for some $p \in \ZZ[q^{\pm 1/2}]$.
Because of the symmetry between~$M$ and~$M'$, the element~$p$ is
invertible, so we conclude that $M'(c) = q^{r/2} \ M(c)$ for some $r \in \ZZ$.
Finally, the fact that $r = 0$ follows from
Proposition~\ref{pr:bar-involution-compatible} since both $M(c)$
and $M'(c)$ are invariant under the bar-involution.
\endproof

\begin{remark}
An important consequence of Theorem~\ref{th:exchange-graphs} is
that the classification of cluster algebras
of finite type achieved in \cite{fz-clust2} applies verbatim to
quantum cluster algebras.
\end{remark}

\begin{remark}
Proposition~\ref{pr:bar-involution-compatible} has the following
important corollary: all cluster variables in $\AA(\mathcal{S})$
are invariant under the bar-involution associated to~$\mathcal{S}$.
A good illustration for this is provided by
Example~\ref{ex:quantum-rank-2}: indeed, the elements given by
\eqref{eq:Ym-Laurent-A2}-\eqref{eq:Ym-Laurent-G2} are obviously
invariant under the bar-involution.
\end{remark}

We conclude this section by exhibiting a family of gradings of the
upper cluster algebras.

\begin{definition}
\label{def:graded-seed}
A \emph{graded quantum seed} is a triple
$(M, \tilde B, \Sigma)$, where
\begin{itemize}
  \item $(M, \tilde B)$ is a quantum seed in $\FF$;
  \item $\Sigma$ is a symmetric integer $m \times m$ matrix such
  that $\tilde B^T \Sigma = 0$.
\end{itemize}
As in Definitions~\ref{def:seed} and \ref{def:quantum-seed},
graded quantum seeds are defined up to a permutation of
the standard basis in $\ZZ^m$ together with the corresponding relabeling of rows and columns of
$\tilde B$ and $\Sigma$.
\end{definition}

We identify $\Sigma$ with the corresponding symmetric bilinear form on $\ZZ^m$.
Then the condition $\tilde B^T \Sigma = 0$ is equivalent to
\begin{equation}
\label{eq:bk-Sigma}
b^j \in \ker \Sigma \quad (j \in \ex),
\end{equation}
where $b^j \in \ZZ^m$ is the $j$th column of $\tilde B$.

The choice of the term ``graded" in
Definition~\ref{def:graded-seed} is justified by the following
construction: every graded quantum seed $(M, \tilde B, \Sigma)$
gives rise to a $\ZZ$-grading on the $\ZZ[q^{\pm 1/2}]$-module $\TT_M$ given by
\begin{equation}
\label{eq:grading-sigma}
{\rm deg}_\Sigma (M(c)) = \Sigma(c,c) \quad (c \in \ZZ^m) \ .
\end{equation}
(Note that this is \emph{not} an algebra grading.)

We will extend quantum seed mutations to graded quantum seeds.
Fix an index $k \in \ex$ and a sign $\varepsilon \in \{\pm 1\}$.
Let $\tilde B'$ be obtained from $\tilde B$ by the mutation in
direction~$k$, and set
\begin{equation}
\label{eq:Sigma-mutated}
 \Sigma' = E_\varepsilon^T \Sigma E_\varepsilon \ ,
\end{equation}
where $E_\varepsilon$ has the same meaning as in \eqref{eq:E-entries}.
Clearly, $\Sigma'$ is symmetric.
The following proposition is an analogue of Proposition~\ref{pr:triple-mutation}
and is proved by the same argument.

\begin{proposition}
\label{pr:Sigma-mutation}
\begin{enumerate}
\item We have $(\tilde B')^T \ \Sigma' = 0$.

\item $\Sigma'$ is independent of the choice of a sign~$\varepsilon$.
\end{enumerate}
\end{proposition}

Proposition~\ref{pr:Sigma-mutation} justifies the following
definition, which extends Definition~\ref{def:quantum-seed-mutation}.

\begin{definition}
\label{def:graded-seed-mutation}
Let $(M, \tilde B, \Sigma)$ be a graded quantum seed, and $k \in \ex$.
We say that the graded quantum seed $(M', \tilde B', \Sigma')$
is obtained from $(M, \tilde B, \Sigma)$
by the \emph{mutation} in direction~$k$,
and write $(M', \tilde B', \Sigma') = \mu_k (M, \tilde B, \Sigma)$
if $(M', \tilde B') = \mu_k (M, \tilde B)$, and
$\Sigma'$ is given by \eqref{eq:Sigma-mutated}.
\end{definition}

Clearly, the mutations of graded quantum seeds are involutive
(cf. Proposition~\ref{pr:quantum-mutation-involutive}).
Therefore, we can define the mutation-equivalence for graded quantum seeds,
and the exchange graph $E(\tilde {\mathcal{S}})$ for a mutation-equivalence class of
graded quantum seeds in the same way as for ordinary quantum seeds above.

\begin{proposition}
\label{pr:grading-rigid}
Let $\tilde {\mathcal{S}}$ be the mutation-equivalence class of
a graded quantum seed $(M, \tilde B, \Sigma)$, and
$\mathcal{S}$ be the mutation-equivalence class of
the underlying quantum seed $(M, \tilde B)$.
\begin{enumerate}
  \item The upper cluster algebra $\mathcal{U}(\mathcal{S})$
  is a graded $\ZZ[q^{\pm 1/2}]$-submodule of $(\TT_M, {\rm deg}_\Sigma)$; furthermore, the
  restriction of the grading ${\rm deg}_\Sigma$ to $\mathcal{U}(\mathcal{S})$
  does not depend on the choice of a
  representative of $\tilde {\mathcal{S}}$.
  \item The forgetful map $(M, \tilde B, \Sigma) \mapsto (M, \tilde B)$
  is a bijection between $\tilde {\mathcal{S}}$ and $\mathcal{S}$,
  i.e., it identifies the exchange graph $E(\tilde {\mathcal{S}})$
  with $E(\mathcal{S})$.
\end{enumerate}
\end{proposition}

\begin{proof}
As in the proof of Proposition~\ref{pr:bar-involution-compatible},
to prove (1) it suffices to show the following: if two graded quantum seeds
$(M, \tilde B, \Sigma)$ and $(M', \tilde B', \Sigma')$
are obtained from each other by a mutation in some direction~$k$, then
$\TT_M \cap \TT_{M'}$ is a graded
$\ZZ[q^{\pm 1/2}]$-submodule of each of
$(\TT_M, {\rm deg}_\Sigma)$ and $(\TT_{M'}, {\rm deg}_{\Sigma'})$,
and the restrictions of both gradings to $\TT_M \cap \TT_{M'}$ are the same.
By the same argument as in the proof of
Proposition~\ref{pr:bar-involution-compatible},
it is enough to show that, for every $c \in \ZZ^m$ with $c_k \geq 0$,
the element $M'(c) \in \TT_M \cap \TT_{M'}$ is homogeneous with
respect to ${\rm deg}_\Sigma$, and
${\rm deg}_{\Sigma}(M'(c)) = \Sigma'(c,c)$.
By \eqref{eq:M'}, $M'(c)$ is a $\ZZ[q^{\pm 1/2}]$-linear
combination of the elements
$M(E_\varepsilon c + \varepsilon p b^k)$; to complete the proof of~(1),
it remains to note that, in view of \eqref{eq:bk-Sigma} and
\eqref{eq:Sigma-mutated}, we have
$$\Sigma(E_\varepsilon c + \varepsilon p b^k, E_\varepsilon c + \varepsilon p b^k)
= \Sigma(E_\varepsilon c, E_\varepsilon c) = \Sigma'(c,c),$$
as required.

\medskip

To prove (2), suppose that $\tilde {\mathcal{S}}$ contains two
graded quantum sets $(M, \tilde B, \Sigma)$ and
$(M, \tilde B, \Sigma')$ with the same underlying quantum seed.
By the already proven part~(1), the two gradings
${\rm deg}_{\Sigma}$ and ${\rm deg}_{\Sigma'}$ agree with each
other on $\mathcal{U}(\mathcal{S})$.
In particular, for every $c \in \ZZ_{\geq 0}^m$, we have
$$\Sigma (c,c) = {\rm deg}_{\Sigma}(M(c)) =
{\rm deg}_{\Sigma'}(M(c)) = \Sigma' (c,c) \ .$$
It follows that $\Sigma = \Sigma'$, and we are done.
\end{proof}

Proposition~\ref{pr:grading-rigid} allows us to include the
bar-involution on $\mathcal{U}(\mathcal{S})$ into a family of more
general ``twisted" bar-involutions defined as follows.
Let $(M, \tilde B, \Sigma)$ be a graded quantum seed.
We associate with $(M,\tilde B, \Sigma)$ the $\ZZ$-linear \emph{twisted bar-involution}
$X \mapsto \overline {X}^{(\Sigma)}$ on $\TT_M$ by the following
formula generalizing \eqref{eq:bar-involution-trivial}:
\begin{equation}
\label{eq:twisted-bar-involution}
\overline {q^{r/2} M(c)}^{(\Sigma)} = q^{-(r+\Sigma(c,c))/2} M(c) \quad (r \in \ZZ, \, c \in \ZZ^m)\ .
\end{equation}
The following proposition generalizes
Proposition~\ref{pr:bar-involution-compatible}.

\begin{proposition}
\label{pr:twisted-bar-involution}
The twisted bar-involution $X \mapsto \overline {X}^{(\Sigma)}$
associated with a graded quantum seed $(M,\tilde B, \Sigma)$
preserves the subalgebra $\mathcal{U}(M, \tilde B)$ of $\TT_M$,
and its restriction to $\mathcal{U}(M, \tilde B)$ depends only
on the mutation-equivalence class of $(M,\tilde B, \Sigma)$.
\end{proposition}

\begin{proof}
Recall the $\ZZ$-grading ${\rm deg}_{\Sigma}$ on $\TT_M$ given by
\eqref{eq:grading-sigma}, and note that
the twisted bar-involution $X \mapsto \overline {X}^{(\Sigma)}$ on
$\TT_M$ can be written as follows:
\begin{equation}
\label{eq:twisted-bar-involution-2}
\overline {X}^{(\Sigma)} = Q^{-1} (\overline {Q(X)}) \ ,
\end{equation}
where~$Q$ is a $\ZZ[q^{\pm 1/2}]$-linear map given by
$Q(X) = q^{d/4} X$ for every homogeneous element $X \in \TT_M$ of
degree~$d$.
By Part~(1) of Proposition~\ref{pr:grading-rigid}, the map~$Q$
preserves the subalgebra $\mathcal{U}(M, \tilde B) \subset \TT_M$,
and its restriction to $\mathcal{U}(M, \tilde B)$ depends only
on the mutation-equivalence class of $(M,\tilde B, \Sigma)$.
Therefore, the same is true for the twisted bar-involution.
\end{proof}

\section{Lower bounds and acyclicity}
\label{sec:lower}

In this section we state and prove quantum analogs of the results
in \cite{bfz-clust3} concerning \emph{lower bounds.}
We retain all the notation and assumptions in
Section~\ref{sec:upper-bound}.
In particular, we assume (without loss of generality) that $L = \ZZ^m$, and the
toric frame $M$ of the ``initial" quantum seed
$(M,\tilde B)$ is given by $M(c) = X^c$ for $c \in L$.
Furthermore, we assume that the initial cluster $\mathbf{X}$ is
the set $\{X_1, \dots, X_n\}$, where $X_j = X^{e_j}$.
By \eqref{eq:quantum-exchange}, for $k \in [1,n]$, the mutation in direction~$k$ replaces
$X_k$ with an element $X'_k$ given by
\begin{equation}
\label{eq:X'k-standard}
X'_k = X^{- e_k + \sum_{b_{ik} > 0} b_{ik} e_i} +
X^{- e_k - \sum_{b_{ik} < 0} b_{ik} e_i} \ .
\end{equation}
It follows that $X'_k$ quasi-commutes with all $X_i$ for $i \neq k$; and each of
the products $X_k X'_k$ and $X'_k X_k$ is the sum of two monomials
in $X_1, \dots, X_m$.
The elements $X'_1, \dots, X'_n$ also
satisfy the following (quasi-)commutation relations.

\begin{proposition}
\label{pr:X'-q-commutation}
Let $j$ and $k$ be two distinct indices from $[1,n]$.
Then $X'_j X'_k - q^{r/2} X'_k X'_j = (q^{s/2} - q^{t/2}) X^e$
for some integers $r,s,t$, and some vector $e \in \ZZ_{\geq 0}^m$.
\end{proposition}

\begin{proof}
Without loss of generality, assume that $b_{jk} \leq 0$.
We abbreviate
$$e'_j = - e_j + \sum_{b_{ij} > 0} b_{ij} e_i, \quad
e'_k = - e_k - \sum_{b_{ik} < 0} b_{ik} e_i,$$
so that \eqref{eq:X'k-standard} can be rewritten as
$$X'_j = X^{e'_j} + X^{e'_j - b^j}, \quad
X'_k = X^{e'_k} + X^{e'_k + b^k},$$
where the vectors $b^j, b^k \in \ZZ^m$ are the $j$th and $k$th columns
of~$\tilde B$.
Using \eqref{eq:multiplication-F} and \eqref{eq:bk-properties}, we
obtain
\begin{eqnarray}
\label{eq:X'-q-commute-explicit}
q^{- \Lambda(e'_j - b^j,e'_k+b^k)/2} X'_j X'_k -
q^{- \Lambda(e'_k+b^k,e'_j - b^j)/2} X'_k X'_j &  \\
\nonumber
= (q^{- \Lambda(b^j,b^k)/2} - q^{- \Lambda(b^k,b^j)/2}) X^{e'_j+e'_k}  \ . &
\end{eqnarray}
If $b_{jk} = 0$ then $\Lambda(b^j,b^k) = 0$ by
\eqref{eq:bk-properties}, and so the right hand side of
\eqref{eq:X'-q-commute-explicit} is equal to $0$; we see that in
this case, $X'_j$ and $X'_k$ quasi-commute.
And if $b_{jk} < 0$ (and so $b_{kj} > 0$) then the vector
$e = e'_j+e'_k$ belongs to $\ZZ_{\geq 0}^m$,
since its $j$th (resp.~$k$th) component is
$-b_{jk} - 1 \geq 0$ (resp. $b_{kj} - 1 \geq 0$).
\end{proof}

Following \cite[Definition~1.10]{bfz-clust3},
we associate with a quantum seed $(M, \tilde B)$ the algebra
\begin{equation}
\label{eq:lower}
\LL(M, \tilde B) =
\ZZ \PP[X_1, X'_1, \dots, X_n, X'_n] \,.
\end{equation}
We refer to $\LL(M, \tilde B)$ as the \emph{lower bound}
associated with $(M, \tilde B)$; this name is justified by
the obvious inclusion $\LL(M,\tilde B) \subset \AA(M,\tilde B)$.

The following definition is an analog of
\cite[Definition~1.15]{bfz-clust3}.

\begin{definition}
\label{def:standard-monomials}
A \emph{standard monomial} in
$X_1, X'_1, \ldots, X_n, X'_n$ is an element of the form
$X_1^{a_1} \cdots X_n^{a_n} (X'_1)^{a'_1} \cdots (X'_n)^{a'_n}$,
where all exponents are nonnegative integers, and $a_k a'_k = 0$ for
$k \in [1,n]$.
\end{definition}

Using the relations between the elements
$X_1, \dots, X_n, X'_1, \dots, X'_n$
described above, it is easy to see that
\begin{equation}
\label{eq:standard-monomials-span}
\text{the standard monomials generate
$\LL(M, \tilde B)$ as a $\ZZ \PP$-module.}
\end{equation}

To state our first result on the lower bounds, we need to recall
the definition of \emph{acyclicity} given in \cite[Definition~1.14]{bfz-clust3}.
We encode the sign pattern of matrix entries of the exchange matrix~$B$
(i.e., the principal part of $\tilde B$) by the
directed graph $\Gamma (B)$ with the vertices $1, \dots, n$ and
the directed edges $(i,j)$ for $b_{ij} > 0$.
We say that $B$ (as well as the corresponding quantum seed)
is \emph{acyclic} if $\Gamma (B)$ has no oriented cycles.
The following result is an analog of
\cite[Theorem~1.16]{bfz-clust3}.

\begin{theorem}
\label{th:standard-monomials-acyclic}
The standard monomials in $X_1, X'_1, \ldots,
X_n, X'_n$ are linearly independent over~$\ZZ \PP$
(that is, they form a $\ZZ \PP$-basis of
$\LL(M,\tilde B)$) if and only if $B$ is acyclic.
\end{theorem}

\begin{proof}
The proof goes along the same lines as that of
\cite[Theorem~1.16]{bfz-clust3}.
The only place where one has to be a little careful is
\cite[Lemma~5.2]{bfz-clust3} which is modified as follows.

\begin{lemma}
Let $u_1, \dots, u_\ell$ and $v_1, \dots, v_\ell$ be some elements of
an associative ring, and let
$i \mapsto i^+$ be a cyclic permutation of~$[1,\ell]$.
For every subset $J\subset [1,\ell]$ such that
$J\cap J^+ = \emptyset$, and for every $i \in [1,\ell]$, we set
$$t_i (J) =
\begin{cases}
u_i & \text {if $i \in J$,} \\
v_i & \text {if $i \in J^+$.}\\
u_i + v_i & \text {otherwise.}
\end{cases}
$$
Then
\begin{equation}
\label{eq:cyclic-identity}
\doublesubscript{\sum}{J\subset [1,\ell]}{J\cap J^+=\emptyset} (-1)^{|J|}
t_1 (J) \cdots t_\ell (J) =
u_1 \cdots u_\ell + v_1 \cdots v_\ell \ .
\end{equation}
\end{lemma}

The proof of \cite[Lemma~5.2]{bfz-clust3} applies verbatim, and
the rest of the proof of \cite[Theorem~1.16]{bfz-clust3}
holds with obvious modifications.
\end{proof}

Our next result is an analog of \cite[Theorem~1.18]{bfz-clust3};
it shows that the acyclicity condition closes the gap
between the upper and lower bounds.

\begin{theorem}
\label{th:acyclic-closing-gaps}
If a quantum seed $(M,\tilde B)$ is acyclic then
$\LL(M, \tilde B) =
\AA(\mathcal {S}) = \mathcal{U}(\mathcal{S})$, where
$\mathcal{S}$ is the mutation-equivalence class of $(M, \tilde B)$.
\end{theorem}

\begin{proof}
The proof of \cite[Theorem~1.18]{bfz-clust3}
extends to the quantum setting, again with some modifications
caused by non-commutativity.
The most non-trivial of these modifications is the following: in
\cite[Lemma~6.7]{bfz-clust3}, we have to replace $P_1$ with an
element $P_{b^1,+}^r$ for an arbitrary positive integer~$r$;
the proof of the modified claim then follows from
Proposition~\ref{pr:center} in the same way
as in Case~1 in the proof of
Proposition~\ref{pr:upper-bound-intersection-2}.
\end{proof}

We conclude this section with an analog of
 \cite[Theorem~1.20]{bfz-clust3},
 which is proved in the same way as its prototype.

\begin{theorem}
\label{th:acyclic-criterion-L=A}
The condition that a quantum seed
$(M, \tilde B)$ is acyclic, is necessary and sufficient for the equality
$\LL(M, \tilde B) = \AA(\mathcal {S})$.
\end{theorem}

\section{Matrix triples associated with Cartan matrices}
\label{sec:Cartan-triples}

In this section we construct a class of matrix triples
$(\Lambda, \tilde B, \Sigma)$ satisfying conditions in
Definitions~\ref{def:seed}, \ref{def:compatible-triple} and
\ref{def:graded-seed}, i.e., giving rise to graded quantum seeds
in the sense of Definition~\ref{def:graded-seed}.
These triples are associated with (generalized) Cartan matrices;
in the case of finite type Cartan matrices, the matrices $\tilde B$
were introduced in \cite[Definition~2.3]{bfz-clust3}.
Our terminology on Cartan matrices and related notions will
basically follow \cite{kac}.

\subsection{Cartan data}
\label{sec:cartan-data}

\begin{definition}
\label{def:cartan-matrix}
A (generalized) \emph{Cartan matrix} is an $r\times r$ integer matrix $A=(a_{ij})$
such that
\begin{itemize}
\item $a_{ii}=2$ for all $i$.

\item $a_{ij} \leq 0$ for $i \neq j$.

\item $a_{ij} = 0$ if and only if $a_{ji} = 0$.

\end{itemize}
\end{definition}

Recall that $A$ is \emph{symmetrizable} if
$d_i a_{ij} = d_j a_{ji}$ for
some positive integers $d_1, \dots, d_r$.
In what follows, we fix a symmetrizable Cartan matrix $A$
and the numbers $d_i$.

\begin{definition}
\label{def:realization}
A \emph{realization} of~$A$ is a triple
$(\hh, \Pi, \Pi^\vee)$, where $\hh$ is a $\CC$-vector space,
and $\Pi = \{\alpha_1, \dots, \alpha_r\} \subset \hh^*$, and
$\Pi^\vee = \{\alpha_1^\vee, \dots, \alpha_r^\vee\} \subset \hh$
are two subsets satisfying the following conditions:
\begin{itemize}
\item both $\Pi$ and $\Pi^\vee$ are linearly independent.

\item $\alpha_j (\alpha_i^\vee) = a_{ij}$ for all $i, j$.

\item $\dim \hh + {\rm rk} A = 2r$.
\end{itemize}
\end{definition}

In what follows, we fix a realization of~$A$; as shown in \cite[Proposition~1.1]{kac},
it is unique up to an isomorphism.
The elements $\alpha_i$ (resp.~$\alpha_i^\vee$) are called \emph{simple roots}
(resp. \emph{simple coroots}) associated to~$A$.

For each $i \in [1,r]$, the \emph{simple reflection} $s_i$ is an
involutive linear transformation of $\hh^*$ acting by
$$s_i (\gamma) = \gamma - \gamma (\alpha_i^\vee) \alpha_i.$$
The Weyl group $W$ is the group generated by all~$s_i$.
We fix a family $\{\omega_1, \dots, \omega_r\} \subset \hh^*$ such
that $\omega_j (\alpha_i^\vee) = \delta_{ij}$; the elements
$\omega_j$ are called \emph{fundamental weights}.
Thus, we have
\begin{equation}
\label{eq:s-omega}
s_i (\omega_j) =
\begin{cases}
\omega_j - \alpha_j & \text{if $i=j$;} \\
\omega_j  & \text{if $i \neq j$.}
\end{cases}
\end{equation}
Note that each $\omega_j$ is defined up to a translation by a
$W$-invariant vector from $\hh^*$.
Note also the following useful property:
\begin{equation}
\label{alpha-omega}
\text{for every $j \in [1,r]$, the vector
$
\sum_{i \in [1,r]} a_{ij} \omega_i - \alpha_j$
is $W$-invariant.}
\end{equation}

As shown in \cite[Chapter~2]{kac}, there exists a $W$-invariant
nondegenerate symmetric bilinear form $(\gamma | \delta)$ on
$\hh^*$ such that
\begin{equation}
\label{eq:scalar-with-alpha}
(\alpha_i | \gamma) = d_i \gamma (\alpha_i^\vee)
\quad (\gamma \in \hh^*).
\end{equation}

\subsection{Double words and associated matrix triples}
\label{sec:double-words}
By a \emph{double word} we will mean a sequence $\ii=(i_1,\ldots,i_m)$ of indices
from $\pm [1,r] = - [1,r] \sqcup [1,r]$.
For every $i \in [1,r]$, we denote
$$\varepsilon (\pm i) = \pm 1, \quad |\pm i| = i.$$
We adopt the convention that $s_{-i}$ is the identity
transformation of $\hh^*$ for $i \in [1,r]$.
For any $a \leq b$ in $[1,m]$, and any sign $\varepsilon$, we set
$$\pi_\varepsilon [a,b] = \pi_\varepsilon^\ii [a,b] =
s_{\varepsilon i_a} \cdots s_{\varepsilon i_b} \ .$$
Iterating \eqref{eq:s-omega}, we obtain the following properties
which will be used many times below:
\begin{eqnarray}
\label{eq:pi-reduction}
&&\text{$\pi_\varepsilon [a,b] \omega_i =
\pi_\varepsilon [a,c] \omega_i$ if $a \leq c \leq b$, and
$\varepsilon i_t \neq  i$ for $c < t \leq b$,}\\
\nonumber
&&\text{$\pi_\varepsilon [a,b] \omega_j =
\pi_\varepsilon [a,b-1] (\omega_j - \alpha_j)$
if $\varepsilon i_b = j$.}
\end{eqnarray}

For $k \in [1,m]$, we denote by $k^+ = k^+_\ii$
the smallest index $\ell$ such that $k < \ell \leq m$ and
$|i_\ell| = |i_k|$; if $|i_k| \neq |i_\ell|$ for $k < \ell \leq m$, then
we set $k^+ = m+1$.
Let $k^- = k^-_\ii$ denote the index $\ell$ such that $\ell^+=k$;
if such an $\ell$ does not exist, we set $k^- = 0$.
We say that an index $k\in [1,m]$ is \emph{$\ii$-exchangeable} if
both $k^-$ and $k^+$ belong to $[1,m]$, and denote by
$\ex = \ex_\ii \subset [1,m]$ the subset of
$\ii$-exchangeable indices.

We will associate to a double word $\ii$ a triple
$(\Lambda (\ii), \tilde B (\ii), \Sigma (\ii))$, where
$\Lambda (\ii)$ and  $\Sigma (\ii)$ are integer $m \times m$
matrices (respectively, skew-symmetric and symmetric), while
$\tilde B (\ii)$ is a rectangular integer matrix with rows
labeled by $[1,m]$ and columns labeled by~$\ex$.

We define the matrix entries of $\Lambda (\ii)$ and  $\Sigma (\ii)$ by
\begin{equation}
\label{eq:lambda-sigma-ii}
\lambda_{k \ell} = \eta_{k, \ell^+} - \eta_{\ell, k^+}, \quad
\sigma_{k \ell} = \eta_{k, \ell^+} + \eta_{\ell, k^+}
\end{equation}
for $k, \ell \in [1,m]$, where
\begin{equation}
\label{eq:eta-kl}
\eta_{k \ell} = \eta_{k \ell}(\ii) =
(\pi_- [\ell, k] \omega_{|i_k|} - \pi_+ [\ell, k] \omega_{|i_k|}
 \ | \ \omega_{|i_\ell|})
\end{equation}
(with the convention that $\eta_{k \ell} = 0$ unless
$1 \leq \ell \leq k \leq m$).
Note that $\eta_{k \ell}$ and so both matrices $\Lambda (\ii)$ and
$\Sigma (\ii)$ are independent of the choice of fundamental
weights.
Indeed, a simple calculation shows that $\eta_{k \ell}$ does not
change if we replace $\omega_{|i_k|}$ by
$\omega_{|i_k|} + \gamma$, and $\omega_{|i_\ell|}$ by
$\omega_{|i_\ell|} + \gamma'$, where both $\gamma$ and $\gamma'$
are $W$-invariant.

Following \cite[Definitions~2.2, 2.3]{bfz-clust3} (which in turn
were based on \cite{z-imrn}), we define the matrix entries
$b_{pk}$ of $\tilde B (\ii)$ for $p \in [1,m]$ and
$k \in \ex$ as follows:
\begin{equation}
\label{eq:tildeB-entries}
b_{pk} = b_{pk}(\ii) =
\begin{cases}
-\varepsilon (i_k) & \text{if $p = k^-$;} \\
-\varepsilon (i_k) a_{|i_p|, |i_k|} &
\text{if $p < k < p^+ < k^+, \, \varepsilon (i_k) = \varepsilon(i_{p^+})$,}\\
 & \text{or $p < k < k^+ < p^+, \, \varepsilon (i_k) = -\varepsilon (i_{k^+})$;} \\
\varepsilon (i_p) a_{|i_p|, |i_k|} &
\text{if $k < p < k^+ < p^+, \, \varepsilon (i_p) = \varepsilon(i_{k^+})$,}\\
& \text{or $k < p < p^+ < k^+, \, \varepsilon (i_p) = -\varepsilon (i_{p^+})$;} \\
\varepsilon (i_p) & \text{if $p = k^+$;} \\
0 & \text{otherwise.}
\end{cases}
\end{equation}
(For technical reasons, the matrix $\tilde B (\ii)$ given by
\eqref{eq:tildeB-entries} differs by sign from the one in
\cite[Definitions~2.2, 2.3]{bfz-clust3}, but this does not affect
the corresponding cluster algebra structure.)

\begin{theorem}
\label{th:compatible-ii}
Suppose that a double word $\ii$
satisfies the following condition:
\begin{eqnarray}
\label{eq:ii-accept-stronger}
&&\text{for every $p \in [1,m]$ with $p^- = 0$, there are no}\\
\nonumber
&&\text{$\ii$-exchangeable indices $k \in [1,p-1]$ with
$a_{|i_p|,|i_k|} < 0$.}
\end{eqnarray}
Then the matrix entries given by {\rm \eqref{eq:lambda-sigma-ii}}
and {\rm \eqref{eq:tildeB-entries}} satisfy
\begin{equation}
\label{eq:b-lambda-sigma}
\sum_{p=1}^m b_{pk} \lambda_{p \ell} = 2 \delta_{k \ell}d_{|i_k|},
\quad \sum_{p=1}^m b_{pk} \sigma_{p \ell} = 0
\end{equation}
for $\ell \in [1,m]$ and $k \in \ex$.
Thus the pair $(\Lambda(\ii),\tilde B(\ii))$ is compatible
in the sense of Definition~\ref{def:compatible-triple},
and the pair $(\tilde B(\ii),\Sigma(\ii))$ satisfies
Definition~\ref{def:graded-seed}.
\end{theorem}

\begin{example}
\label{ex:ii-standard-sl3}
Let
$$A =
  \begin{pmatrix}
2 &  -1 \\
-1 &  2
  \end{pmatrix},
$$
be the Cartan matrix of type $A_2$, with $d_1 = d_2 = 1$.
Taking
$$\ii = (1,2,1,2,1,-1,-2,-1),$$
it is easy to check that the corresponding
matrices $\tilde B(\ii)$ and $\Lambda(\ii)$ are those in
Example~\ref{ex:sl3-quantized-matrices}.
The first equality in \eqref{eq:b-lambda-sigma} was shown there.
As for $\Sigma(\ii)$, it is a symmetric matrix whose entries on
and below the main diagonal are equal to those of $\Lambda(\ii)$.
The last equality in \eqref{eq:b-lambda-sigma} can be seen by a
direct inspection.
\end{example}

\proof
We will use \eqref{eq:tildeB-entries} to define
$b_{pk}$ for all $k, p \in [1,m]$ (with $k$ not necessarily
$\ii$-exchangeable).
In view of \eqref{eq:lambda-sigma-ii}, to verify
\eqref{eq:b-lambda-sigma} it suffices to show the following.

\begin{lemma}
\label{lem:b-eta}
For an arbitrary double word $\ii$, we have
\begin{equation}
\label{eq:b-eta1}
\sum_{p=1}^m b_{pk} \eta_{p \ell} = \delta_{k^+, \ell}\ d_{|i_k|}
\end{equation}
for all $k, \ell \in [1,m]$ such that $k^+ \leq m$.
If $\ii$ satisfies {\rm \eqref{eq:ii-accept-stronger}} then we
also have
\begin{equation}
\label{eq:b-eta2}
\sum_{p=1}^m b_{pk} \eta_{\ell, p^+} = - \delta_{k \ell} d_{|i_k|}
\end{equation}
for all $\ell \in [1,m]$ and $k \in \ex$.
\end{lemma}

The rest of this section is dedicated to the proof of
Lemma~\ref{lem:b-eta}.
First, we get \eqref{eq:b-eta2} out of the way by showing that it
follows from \eqref{eq:b-eta1}.
To see this, consider the \emph{opposite} double word
$\ii^\circ = (i_m, \dots, i_1)$.
We abbreviate $k^\circ = m+1-k$, so that $\ii^\circ$ can be
written as $\ii^\circ = (i_{1^\circ}, \dots, i_{m^\circ})$.
Examining \eqref{eq:eta-kl} and \eqref{eq:tildeB-entries}, we obtain
\begin{align}
\label{eq:eta-b-opposite}
  \eta_{k \ell}(\ii) & = \eta_{\ell^\circ, k^\circ}(\ii^\circ)
  \quad (k, \ell \in [1,m]), \\
\nonumber
  b_{p k}(\ii) & = - b_{{p^+}^\circ, {k^+}^\circ}(\ii^\circ) \quad
(k^+, p^+ \in [1,m]).
\end{align}

Turning to \eqref{eq:b-eta2}, we note that the summation there
can be restricted to the values of~$p$ such that $p^+ \leq m$
(because $\eta_{\ell, p^+} = 0$ unless $p^+ \leq \ell$).
Substituting the expressions given by \eqref{eq:eta-b-opposite}
into \eqref{eq:b-eta2}, we obtain
\begin{equation}
\label{eq:b-eta2-opposite}
\sum_{p=1}^m b_{pk} \eta_{\ell, p^+} =
- \sum_{p^+ \leq m} b_{{p^+}^\circ, {k^+}^\circ}(\ii^\circ)
\eta_{{p^+}^\circ, \ell^\circ}(\ii^\circ).
\end{equation}
Comparing this with the counterpart of \eqref{eq:b-eta1} for the double word
$\ii^\circ$, we see that it remains to show the following:
$$\sum_{{(p^\circ)}^+_{\ii^\circ} = m+1} b_{p^\circ, {k^+}^\circ}(\ii^\circ)
\eta_{p^\circ, \ell^\circ}(\ii^\circ) = 0,$$
whenever $k$ is $\ii$-exchangeable.
To complete the proof of \eqref{eq:b-eta2}, it remains to observe
that the condition \eqref{eq:ii-accept-stronger} guarantees that
$b_{p^\circ, {k^+}^\circ}(\ii^\circ) = 0$ for all $p$ such that
${(p^\circ)}^+_{\ii^\circ} = m+1$ (which is equivalent to $p^- = 0$).

We now concentrate on the proof of \eqref{eq:b-eta1}.
We will need to consider several cases of the relative position of~$k$ and~$\ell$.
As a warm-up, we note that $b_{pk} = 0$ for $p > k^+$, and
$\eta_{p \ell} = 0$ for $p < \ell$; therefore, the sum in
\eqref{eq:b-eta1} is equal to~$0$ if $\ell > k^+$.
For $\ell = k^+$, the sum in question reduces to just one term
with $p = \ell = k^+$; using \eqref{eq:eta-kl}, \eqref{eq:tildeB-entries},
and \eqref{eq:s-omega}-\eqref{eq:scalar-with-alpha},
we see that this term is equal to
\begin{align*}
b_{pk} \eta_{p \ell} &=  \varepsilon (i_p)(s_{-i_p}\omega_{|i_p|} - s_{i_p}
\omega_{|i_p|}\ | \ \omega_{|i_p|})
= (\omega_{|i_k|} - s_{|i_k|} \omega_{|i_k|}\ | \ \omega_{|i_k|})\\
& = (\alpha_{|i_k|}\ | \ \omega_{|i_k|}) = d_{|i_k|},\\
\end{align*}
in accordance with \eqref{eq:b-eta1}.

For the rest of the proof, we assume that $\ell < k^+$, and
(for typographical reasons) abbreviate $|i_k| = j$ and
$|i_\ell| = h$.
To show that the sum in \eqref{eq:b-eta1} is equal to~$0$,
we compute, for every $i \in [1,r]$, the contribution
to this sum from the values of~$p$ such that $|i_p| = i$.
We denote this contribution by~$S_i = S_i(k,\ell;\ii)$.

\begin{lemma}
\label{lem:contributions}
We have
\begin{equation}
\label{eq:Sj}
S_j =
\begin{cases}
(\omega_j - \pi_{\varepsilon(i_{k^+})}[\ell,k^+]\omega_j \ | \ \omega_h)
 & \text{if $k < \ell < k^+$;}\\[.05in]
(\pi_{\varepsilon(i_{k})}[\ell,k](\omega_j - \alpha_j)&\\
- \pi_{\varepsilon(i_{k^+})}[\ell,k^+]\omega_j \ | \ \omega_h)
 & \text{if $\ell \leq k,
 \, \varepsilon (i_k) = \varepsilon(i_{k^+})$;} \\[.05in]
(\pi_{\varepsilon(i_{k})}[\ell,k]\omega_j&\\
- \pi_{\varepsilon(i_{k^+})}[\ell,k^+]\omega_j \ | \ \omega_h)
 & \text{if $k^- < \ell \leq k,
 \, \varepsilon (i_k) = -\varepsilon(i_{k^+})$;}\\[.05in]
(\pi_{\varepsilon(i_{k})}[\ell,k](2\omega_j - \alpha_j)&\\
 -\pi_{\varepsilon(i_{k^+})}[\ell,k^+](2\omega_j - \alpha_j) | \omega_h)
 & \text{if $\ell \leq k^-,
 \ \varepsilon (i_k) = -\varepsilon(i_{k^+})$,}
\end{cases}
\end{equation}
and, for $i \neq j$,
\begin{equation}
\label{eq:Si-neq-j}
S_i =
\begin{cases}
a_{ij}(\omega_i - \pi_{\varepsilon(i_{k^+})}[\ell,k^+]\omega_i \ | \ \omega_h)
 & \text{if $k < \ell < k^+$;}\\[.05in]
a_{ij}(\pi_{\varepsilon(i_{k})}[\ell,k]\omega_i -
\pi_{\varepsilon(i_{k^+})}[\ell,k^+]\omega_i \ | \ \omega_h)
 & \text{if $\ell \leq k$.}
\end{cases}
\end{equation}
\end{lemma}

\begin{proof}
By \eqref{eq:tildeB-entries}, the only possible values of~$p$ contributing
to~$S_j$ are $p = k^+$ and $p = k^-$ (the latter value appears
only when $\ell \leq k^-$).
Let us do the last case in \eqref{eq:Sj} (the other cases are similar):
$\ell \leq k^-, \, \varepsilon (i_k) = -\varepsilon(i_{k^+}) = \varepsilon$.
Applying \eqref{eq:tildeB-entries} and \eqref{eq:eta-kl}, and using \eqref{eq:pi-reduction}, we get
\begin{align*}
b_{k^+,k} \eta_{k^+,\ell} &=
(\pi_{\varepsilon}[\ell,k^+]\omega_j
- \pi_{-\varepsilon}[\ell,k^+]\omega_j \ | \ \omega_h)\\
&= (\pi_{\varepsilon}[\ell,k]\omega_j
- \pi_{-\varepsilon}[\ell,k^+]\omega_j \ | \ \omega_h),
\end{align*}
and
\begin{align*}
b_{k^-,k} \eta_{k^-,\ell} &=
(\pi_{\varepsilon}[\ell,k^-]\omega_j
- \pi_{-\varepsilon}[\ell,k^-]\omega_j \ | \ \omega_h)\\
&= (\pi_{\varepsilon}[\ell,k](\omega_j - \alpha_j)
- \pi_{-\varepsilon}[\ell,k^+](\omega_j - \alpha_j) \ | \ \omega_h),
\end{align*}
which implies our claim.

Turning to \eqref{eq:Si-neq-j}, we will also consider only the
latter case $\ell \leq k$, the former one being similar and simpler.
The indices~$p$ with $|i_p| = i$, which may have a non-zero contribution to~$S_i$,
fall into the following types:
\begin{enumerate}
\item[{\sl Type 1:}]
$\ell \leq p < k < k^+ < p^+, \, \varepsilon (i_k) = -\varepsilon(i_{k^+})$,
or $\ell \leq p < k < p^+ < k^+, \, \varepsilon (i_k) = \varepsilon(i_{p^+})$.
Using \eqref{eq:eta-kl}, \eqref{eq:tildeB-entries}, and
\eqref{eq:pi-reduction}, we see that
the corresponding contribution to $S_i$ is given by
\begin{equation}
\label{eq:b-eta-Si-1}
b_{pk} \eta_{p \ell}=
a_{ij}(\pi_{\varepsilon(i_{k})}[\ell,k]\omega_i -
\pi_{-\varepsilon(i_{k})}[\ell,k]\omega_i \ | \ \omega_h) \ .
\end{equation}

\item[{\sl Type 2:}]
$k < p < p^+ < k^+, \, \varepsilon (i_p) = -\varepsilon(i_{p^+})$,
or $k < p < k^+ < p^+, \, \varepsilon (i_p) = \varepsilon(i_{k^+})$.
The corresponding contribution to $S_i$ is given by
\begin{equation}
\label{eq:b-eta-Si-2}
b_{pk} \eta_{p \ell}=
a_{ij}(\pi_{-\varepsilon(i_{p})}[\ell,p]\omega_i -
\pi_{\varepsilon(i_{p})}[\ell,p]\omega_i \ | \ \omega_h) \ .
\end{equation}
\end{enumerate}
Note that there is at most one index of type~$1$, but there could
be several indices of type~$2$.
We need to show that all the contributions \eqref{eq:b-eta-Si-1}
and \eqref{eq:b-eta-Si-2} add up to
\begin{equation}
\label{eq:Si-last-case}
S_i = a_{ij}(\pi_{\varepsilon(i_{k})}[\ell,k]\omega_i -
\pi_{\varepsilon(i_{k^+})}[\ell,k^+]\omega_i \ | \ \omega_h).
\end{equation}

First suppose that there are no indices $p$ with $|i_p| = i$ between~$k$
and~$k^+$; in particular, there are no indices~$p$ of type~$2$.
In view of \eqref{eq:pi-reduction}, the sum in
\eqref{eq:Si-last-case} can be rewritten as
$$a_{ij}(\pi_{\varepsilon(i_{k})}[\ell,k]\omega_i -
\pi_{\varepsilon(i_{k^+})}[\ell,k]\omega_i \ | \ \omega_h).$$
This expression is easily seen to vanish unless
$\varepsilon (i_k) = -\varepsilon(i_{k^+})$, and there exists a
(unique) index~$p$ of type~$1$; furthermore, in the latter case,
it agrees with \eqref{eq:b-eta-Si-1}.

Next consider the case when there are some indices $p$ with $|i_p| = i$ between~$k$
and~$k^+$, but none of them are of type~$2$.
In other words, all these values of~$p$ have the same sign, say $\varepsilon$,
of~$i_p$, and we also have $\varepsilon(i_{k^+}) = -\varepsilon$.
In this case, the sum in
\eqref{eq:Si-last-case} can be rewritten as
$$a_{ij}(\pi_{\varepsilon(i_{k})}[\ell,k]\omega_i -
\pi_{-\varepsilon}[\ell,k]\omega_i \ | \ \omega_h).$$
Again, this expression vanishes unless
$\varepsilon (i_k) = \varepsilon$, and there exists a
(unique) index~$p$ of type~$1$; and again, in the latter case,
it agrees with \eqref{eq:b-eta-Si-1}.

It remains to treat the case when there are some indices~$p$ of type~$2$.
Let $p(1) < \cdots < p(t)$ be all such indices.
By the definition, we have
$\varepsilon (i_{p(s)}) = -\varepsilon(i_{p(s+1)})$ for
$s = 1, \dots, t-1$, and
$\varepsilon (i_{p(t)}) = \varepsilon(i_{k^+})$.
Furthermore, \eqref{eq:pi-reduction} yields
$\pi_{-\varepsilon(i_{p(s+1)})}[\ell,p(s+1)]\omega_i =
\pi_{\varepsilon(i_{p(s)})}[\ell,p(s)]\omega_i$
for $s = 1, \dots, t-1$.
This shows that the sum of all expressions \eqref{eq:b-eta-Si-2}
allows telescoping, and so is equal to
\begin{equation}
\label{eq:Si-telescoping}
a_{ij}(\pi_{-\varepsilon(i_{p(1)})}[\ell,k]\omega_i -
\pi_{\varepsilon(i_{k^+})}[\ell,k^+]\omega_i \ | \ \omega_h).
\end{equation}
An easy inspection shows that \eqref{eq:Si-telescoping} agrees with
\eqref{eq:Si-last-case} if there are no indices~$p$ of type~$1$.
In the latter case, we must have
$\varepsilon (i_{k}) = \varepsilon(i_{p(1)})$, and so the sum of
expressions in \eqref{eq:Si-telescoping} and \eqref{eq:b-eta-Si-1}
is equal to that in \eqref{eq:Si-last-case}, as desired.
This completes the proof of Lemma~\ref{lem:contributions}.
\end{proof}

To finish the proof of \eqref{eq:b-eta1}, we need to show
that
$$S := S_j + \sum_{i \neq j} S_i = 0$$
in all the cases in Lemma~\ref{lem:contributions}.
Combining \eqref{eq:Sj} and \eqref{eq:Si-neq-j} with \eqref{alpha-omega}, we get
\begin{equation}
\label{eq:S-total}
S =
\begin{cases}
(\alpha_j - \omega_j&\\
- \pi_{\varepsilon(i_{k^+})}[\ell,k^+](\alpha_j-\omega_j) \ | \ \omega_h)
 & \text{if $k < \ell < k^+$,}\\[.05in]
 (\pi_{\varepsilon(i_{k})}[\ell,k](-\omega_j)&\\
- \pi_{\varepsilon(i_{k^+})}[\ell,k^+](\alpha_j-\omega_j \ | \ \omega_h)
 & \text{if $\ell \leq k,
 \, \varepsilon (i_k) = \varepsilon(i_{k^+})$;} \\[.05in]
(\pi_{\varepsilon(i_{k})}[\ell,k](\alpha_j-\omega_j)&\\
- \pi_{\varepsilon(i_{k^+})}[\ell,k^+](\alpha_j-\omega_j) \ | \ \omega_h)
 & \text{if $k^- < \ell \leq k,
 \, \varepsilon (i_k) = -\varepsilon(i_{k^+})$;}\\[.05in]
0 & \text{if $\ell \leq k^-,
 \, \varepsilon (i_k) = -\varepsilon(i_{k^+})$.}
\end{cases}
\end{equation}
It remains to show that $S = 0$ in each of the first three cases in \eqref{eq:S-total}.
In case~$1$, we have
$\pi_{\varepsilon(i_{k^+})}[\ell,k^+](\alpha_j-\omega_j) = -
\omega_j$, and so $S = (\alpha_j \ | \ \omega_h) = 0$.
In case~$2$ (resp.~$3$), we have
$\pi_{\varepsilon(i_{k^+})}[\ell,k^+](\alpha_j-\omega_j) =
\pi_{\varepsilon(i_{k})}[\ell,k](-\omega_j)$
(resp. $\pi_{\varepsilon(i_{k})}[\ell,k](\alpha_j-\omega_j) = - \omega_j
= \pi_{\varepsilon(i_{k^+})}[\ell,k^+](\alpha_j-\omega_j)$),
which again yields $S = 0$.
This completes the proof of \eqref{eq:b-eta1} and hence those of
Lemma~\ref{lem:b-eta} and Theorem~\ref{th:compatible-ii}.
\endproof

\begin{remark}
\label{rem:condition-on-ii}
Inspecting the above proof, we see that the condition
\eqref{eq:ii-accept-stronger} was used only to ensure
that $b_{p^\circ, {k^+}^\circ}(\ii^\circ) = 0$ for all
$\ii$-exchangeable indices~$k$ and all $p$ with $p^- = 0$.
It follows that \eqref{eq:ii-accept-stronger}
can be replaced, for instance, by the following weaker restriction:
\begin{eqnarray}
\label{eqn:ii-accept}
&&\text{For every $p \in [1,m]$ and $j \in [1,r]$ such that
$p^- = 0$, $a_{|i_p|,j} < 0$,}\\
\nonumber
&&\text{and $\{k \in [1,p-1]: |i_k| = j\} = \{k_1 < \cdots < k_t\}$ with $t \geq 2$,}\\
\nonumber
&&\text{we have $\varepsilon(i_{k_2}) = \cdots = \varepsilon(i_{k_t})$;
if $k_t$ is $\ii$-exchangeable
then also}\\
\nonumber
&&\text{$\varepsilon(i_{k_t}) = - \varepsilon(i_p)$.}
\end{eqnarray}
However, the simpler condition \eqref{eq:ii-accept-stronger} is good
enough for our applications.
For instance, it is satisfied whenever the first $r$ terms of $\ii$
are $\pm 1, \dots, \pm r$ arranged in any order; this covers the
class of double words $\ii$ considered in \cite[Section~2]{bfz-clust3}
and in Section~\ref{sec:quantum double Bruhat cells} below.
\end{remark}

\begin{remark}
Because of the fundamental role played by the matrix~$\tilde B$ in
the theory of cluster algebras, it would be desirable to find an
alternative expression to \eqref{eq:tildeB-entries} involving
fewer special cases.
One such expression was given in \cite[Remark~2.4]{bfz-clust3}.
Here we present another expression that seems to be more manageable.
Namely we claim that, for $p \in [1,m]$ and $k \in \ex$,
\eqref{eq:tildeB-entries} is equivalent to
\begin{equation}
\label{eq:b-thru-s}
b_{pk} = s_{pk} - s_{p,k^+} - s_{p^+,k} + s_{p^+,k^+},
\end{equation}
where
\begin{equation}
\label{eq:s}
s_{pk} = \frac
{{\rm sgn}(p-k) (\varepsilon(i_p) + \varepsilon(i_k))}{4} \,
a_{|i_p|,|i_k|},
\end{equation}
and we use the following convention: if $p^+ = m+1$ then the last
two terms in \eqref{eq:b-thru-s} are given by \eqref{eq:s} with
$i_{m+1} = \pm i_p$ (the choice of a sign does not matter).
The proof of \eqref{eq:b-thru-s} is straightforward, and we leave
it to the reader.
\end{remark}

\section{Preliminaries on quantum groups}
\label{sec:quantum groups}

\subsection{Quantized enveloping algebras}
Our standard reference in this section will be \cite{brown-goodearl}.
We start by recalling the definition of the quantized enveloping algebra
associated with a symmetrizable (generalized) Cartan matrix $A =(a_{ij})$.
We fix a realization $(\hh, \Pi, \Pi^\vee)$ of $A$  as in Definition~\ref{def:realization}.
Let $(\gamma | \delta )$ be the inner product on $\hh^*$ defined by (\ref{eq:scalar-with-alpha}).
Define the weight lattice $P$ by
$$\text{$P = \{\lambda\in \hh^*: \lambda(\alpha_i^\vee)\in \ZZ$ for all $i\in [1,r]\}$} \ .$$
The {\it quantized enveloping algebra} $U$ is a $\QQ(q)$-algebra generated
by the elements $E_i$ and $F_i$ for $i \in [1,r]$, and $K_\lambda$ for $\lambda\in P$,
subject to the following relations:
$$K_\lambda K_\mu = K_{\lambda+\mu}, \,\, K_0 = 1$$
for $\lambda, \mu \in P$;
$$K_\lambda E_i =q^{(\alpha_i|\lambda)} E_i K_\lambda,
\,\, K_\lambda F_i =q^{-(\alpha_i|\lambda)} F_i K_\lambda$$
for $i \in [1,r]$ and $\lambda\in P$;
$$E_i F_j - F_j E_i=\delta_{ij}\frac{K_{\alpha_i}- K_{-\alpha_i}}{q^{d_i}-q^{-d_i}}$$
for $i,j \in [1,r]$;
and the {\it quantum Serre relations}
$$\sum_{p=0}^{1-a_{ij}} (-1)^p 
E_i^{[1-a_{ij}-p;i]} E_j E_i^{[p;i]} = 0,$$
$$\sum_{p=0}^{1-a_{ij}} (-1)^p 
F_i^{[1-a_{ij}-p;i]} F_j F_i^{[p;i]} = 0$$
for $i \neq j$, where 
the notation $X^{[p;i]}$ stands for the \emph{divided power}
\begin{equation}
\label{eq:divided-power}
X^{[p;i]} = \frac{X^p}{[1]_i \cdots [p]_i}, \quad
[k]_i = \frac{q^{kd_i}-q^{-kd_i}}{q^{d_i}-q^{-d_i}} \ .
\end{equation}

The algebra $U$ is a $q$-deformation of the universal enveloping algebra of
the Kac-Moody algebra~$\gg$ associated to~$A$, so it is commonly
denoted by $U = U_q(\gg)$.
It has a natural structure of a bialgebra with the comultiplication $\Delta:U\to U\otimes U$
and the counit homomorphism  $\varepsilon:U\to \QQ(q)$
given by
\begin{equation}
\label{eq:coproduct}
\Delta(E_i)=E_i\otimes 1+K_{\alpha_i}\otimes E_i, \,
\Delta(F_i)=F_i\otimes K_{-\alpha_i}+ 1\otimes F_i, \, \Delta(K_\lambda)=
K_\lambda\otimes K_\lambda \ ,
\end{equation}
\begin{equation}
\label{eq:counit}
\varepsilon(E_i)=\varepsilon(F_i)=0, \quad \varepsilon(K_\lambda)=1\ .
\end{equation}
In fact, $U$ is a Hopf algebra with the antipode antihomomorphism
$S: U \to U$ given by
$$S(E_i) = -K_{-\alpha_i} E_i, \,\, S(F_i) = -F_i K_{\alpha_i}, \,\,
S(K_\lambda) = K_{-\lambda},$$
but we will not need this structure.

Let $U^-$ (resp.~$U^0$; $U^+$) be the $\QQ(q)$-subalgebra of~$U$ generated by
$F_1, \dots, F_r$ (resp. by~$K_\lambda \, (\lambda\in P)$; by $E_1, \dots, E_r$).
It is well-known that $U=U^-\cdot U^0\cdot U^+$ (more precisely,
the multiplication map induces an isomorphism $U^-\otimes U^0\otimes U^+ \to U$).

The algebra $U$ is graded by the root lattice $Q$:
\begin{equation}
\label{eq:U-grading}
\text{$U=\bigoplus_{\alpha\in Q} U_\alpha, \quad
U_\alpha = \{u \in U: K_\lambda u K_{-\lambda} = q^{(\lambda\,|\,\alpha)}\cdot u$
for $\lambda\in P\}$.}
\end{equation}
Thus, we have
$${\rm deg} E_i = \alpha_i, \quad {\rm deg} F_i = - \alpha_i, \quad
{\rm deg} K_\lambda = 0 \ .$$

\subsection{The quantized coordinate ring of~$G$}
\label{sec:OqG}
Our next target is the quantized coordinate ring $\OO_q(G)$ (also known as the
\emph{quantum group}) of the group~$G$ associated to the Cartan matrix~$A$.
Since most of the literature on quantum groups deals only with the case
when $A$ is of finite type, we will also restrict our attention to this case
(even though we have little doubt that all the results extend to Kac-Moody groups).
That is, from now on we assume that $A$ is of finite type, i.e.,
it corresponds to a semisimple Lie algebra~$\gg$.
Let $G$ be the simply-connected semisimple group with the Lie algebra~$\gg$.
Following \cite[Chapter I.8]{brown-goodearl}, the
\emph{quantized coordinate algebra} $\OO_q(G)$ can be defined as follows.

First note that $U^*= {\rm Hom}_{\QQ(q)}(U,\QQ(q))$
has a natural algebra structure:
for $f, g\in U^*$, the product $fg$ is defined by
\begin{equation}
\label{eq:coalgebra-product}
fg(u)=(f\otimes g)(\Delta(u)) = \sum f(u_1) g(u_2)
\end{equation}
for all $u\in U$, where we use the Sweedler summation notation
$\Delta(u) = \sum u_1 \otimes u_2$ (cf.~e.g.,
\cite[Section I.9.2]{brown-goodearl}).
The algebra~$U^*$ has the standard $U-U$-bimodule structure
given by
$$(Y \bullet f \bullet X) (u) = f(X u Y)$$
for $f \in U^*$ and $u, X, Y \in U$.
In view of \eqref{eq:coalgebra-product}, we have
\begin{equation}
\label{eq:translation-product}
Y \bullet (fg) \bullet X =
\sum (Y_1 \bullet f \bullet X_1)(Y_2 \bullet g \bullet X_2) \ .
\end{equation}

Let $U^\circ$ be the \emph{Hopf dual} of~$U$ defined by
$$\text{$U^\circ = \{f \in U^* : f(I) = 0$ for some ideal
$I \subset U$ of finite codimension\}.}$$
Then $U^\circ$ is a subalgebra and a $U-U$-sub-bimodule of~$U^*$.

Slightly modifying the definition in \cite[Section I.8.6]{brown-goodearl}, for every
$\gamma, \delta \in P$, we set
\begin{equation}
\label{eq:PP-grading}
\text{$U^\circ_{\gamma, \delta}=
\{f \in U^\circ : K_\mu \bullet f \bullet K_\lambda =
q^{(\lambda | \gamma) + (\mu | \delta)} f$ for $\lambda, \mu \in P\}$.}
\end{equation}
Finally, we define $\OO_q(G)$ as the $P\times P$-graded subalgebra
of $U^\circ$ given by
$$\OO_q(G)=\bigoplus_{\gamma, \delta\in P}  U^\circ_{\gamma, \delta}$$
(from now on, we will denote the homogeneous components of $\OO_q(G)$ by
$\OO_q(G)_{\gamma, \delta}$ instead of $U^\circ_{\gamma,\delta}$).

It is well-known (see e.g., \cite[Theorem I.8.9]{brown-goodearl})
that $\OO_q(G)$ is a domain.

The algebra $\OO_q(G)$ is a $U-U$-sub-bimodule of~$U^\circ$:
according to \cite[Lemma~I.8.7]{brown-goodearl}, we have
$$\text{$Y \bullet \OO_q(G)_{\gamma, \delta} \bullet X \subset
\OO_q(G)_{\gamma - \alpha, \delta+\beta}$
for $X \in U_\alpha, \,\, Y \in U_\beta$.}$$

We now give a more explicit description of $\OO_q(G)$.
Let
$$\text{$P^+ = \{\lambda\in P:\lambda(\alpha_i^\vee) \geq 0$ for all $i\in [1,r]\}$}$$
be the semigroup of dominant weights.
Thus, $P^+$ is a free additive semigroup generated by fundamental
weights $\omega_1, \dots, \omega_r$.
(Since~$A$ is of finite type, the setup in
Section~\ref{sec:cartan-data} simplifies so that simple coroots
(resp.~simple roots) form a basis in $\hh$ (resp.~$\hh^*$), and
the fundamental weights are uniquely determined by the condition
$\omega_j (\alpha_i^\vee) = \delta_{ij}$.)
To every dominant weight $\lambda \in P^+$ we associate an element
$\Delta^\lambda \in U^*$ given by
\begin{equation}
\label{eq:highest-vector}
\Delta^\lambda(FK_\mu E)=\varepsilon(F) q^{(\lambda | \mu)}\varepsilon(E)
\end{equation}
for $F\in U^-$, $E\in U^+$ and $\mu\in P$.
Let $\mathcal{E}_\lambda = U \bullet \Delta^\lambda \bullet U$ be the
$U-U$-sub-bimodule of~$U^*$ generated by $\Delta^\lambda$.
The following presentation of $\OO_q(G)$ was essentially given in
\cite[Section~I.7]{brown-goodearl}.

\begin{proposition}
\label{pr:peter-weyl}
Each element $\Delta^\lambda$ belongs to $\OO_q(G)_{\lambda, \lambda}$,
each subspace $\mathcal{E}_\lambda$ is a finite-dimensional simple $U-U$-bimodule,
and $\OO_q(G)$ has the direct sum decomposition
$$\OO_q(G)= \bigoplus_{\lambda \in P^+} \mathcal{E}_\lambda \ .$$
\end{proposition}


The reason for our choice of the $P \times P$-grading in
$\OO_q(G)$ is the following: we can view $\OO_q(G)$ as a
$U \times U$-module via
$$(X,Y) f = Y \bullet f \bullet X^T,$$
where $X \mapsto X^T$ is the transpose antiautomorphism of the
$\QQ(q)$-algebra~$U$ given by
$$E_i^T = F_i, \quad F_i^T = E_i, \quad K_\lambda^T = K_\lambda \ .$$
The specialization $q = 1$ transforms $\OO_q(G)$ into a $\gg \times \gg$-module,
and $\OO_q(G)_{\gamma, \delta}$ becomes the weight subspace of
weight $(\gamma, \delta)$ under this action.
In particular, under the specialization $q = 1$, the space
$\mathcal{E}_\lambda$ becomes a simple $\gg \times \gg$-module
generated  by the highest vector $\Delta^\lambda$ of weight
$(\lambda, \lambda)$.

Comparing \eqref{eq:PP-grading} with \eqref{eq:U-grading}, we
obtain the following useful property:
\begin{equation}
\label{eq:UO-pairing}
\text{If the pairing $\OO_q(G)_{\gamma, \delta} \times U_\alpha
\to \QQ(q)$ is non-zero then $\alpha = \gamma - \delta$.}
\end{equation}

\subsection{Quantum double Bruhat cells}
For each $i \in [1,r]$, we adopt the notational convention
$$E_{-i} = F_i, \quad s_{-i} = 1$$
(the latter was already used in Section~\ref{sec:double-words}).
For $i \in \pm [1,r] = -[1,r] \sqcup [1,r]$, we denote by $U_i$ the subalgebra
of~$U$ generated by $U^0$ and $E_i$.
For every double word $\ii=(i_1,\ldots,i_m)$
(i.e., a word in the alphabet $\pm [1,r]$), we set
$$U_\ii = U_{i_1} \cdots U_{i_m} \subset U \ .$$
Denote
$$J_\ii:=\{f\in \OO_q(G): f(U_\ii) =  0\},$$
i.e., $J_\ii$ is the orthogonal complement of $U_\ii$ in $\OO_q(G)$.

Clearly, each $U_\ii$ satisfies $\Delta(U_\ii) \subset U_\ii\otimes
U_\ii$, hence $J_\ii$ is a two-sided ideal in $\OO_q(G)$.
In fact, $J_\ii$ is prime, i.e., $\OO_q(G)/J_\ii$ is a domain (see, e.g.,
\cite[Corollary~10.1.10]{joseph}).

Recall that a \emph{reduced word} for $(u,v) \in W\times W$
is a shortest possible double word $\ii=(i_1,\ldots,i_m)$ such that
$$s_{-i_1} \cdots s_{-i_m}=u, \quad s_{i_1} \cdots s_{i_m}=v \ ;$$
thus, $m = \ell(u)+\ell(v)$, where $\ell: W \to \ZZ_{\ge 0}$
is the length function on~$W$.

\begin{proposition}
\label{pr:U-uv}
If $\ii$ and $\ii'$ are reduced words for the same element
$(u,v)\in W\times W$, then $U_\ii=U_{\ii'}$.
\end{proposition}

\begin{proof}
By the well known Tits' lemma, it suffices to check the
statement in the following two special cases:
\begin{enumerate}
  \item $\ii = (i,j,i, \ldots), \, \ii' = (j,i,j, \ldots)$, where
  $i, j \in [1,r]$, and the length of each of $\ii$ and $\ii'$ is
  equal to the order of $s_i s_j$ in $W$;
  \item $\ii = (i,-j), \, \ii' = (-j,i)$, where $i, j \in [1,r]$.
\end{enumerate}
Case~(1) is treated in \cite{Lusztig-problems}, while
Case~(2) follows easily from the commutation relation between
$E_i$ and $F_j$ in $U$.
\end{proof}

In view of Proposition~\ref{pr:U-uv}, for every $u,v \in W$,
we set $U_{u,v} = U_\ii$, and $J_{u,v} = J_\ii$, where
$\ii$ is any reduced word for $(u,v)$.
The algebra $\OO_q(G)/J_{u,v}$ has the following geometric meaning.
Let $H$ be the maximal torus in~$G$ with Lie algebra~$\hh$, and
let~$B$ (resp.~$B_-$) be the Borel subgroup in~$G$ generated
by~$H$ and the root subgroups corresponding to simple roots
$\alpha_1, \dots, \alpha_r$ (resp.~$-\alpha_1, \dots, -\alpha_r$).
Recall that the Weyl group~$W$ is naturally identified with
${\rm Norm}_G(H)/H$.
For $u, v \in W$, let $G^{u,v}$ denote the double Bruhat cell
$BuB \cap B_- v B_-$ in $G$ (for their properties see
\cite{fz-double}).
Let $\overline {G^{u,v}}$ denote the Zariski closure of $G^{u,v}$ in $G$.
As shown in \cite{dc-pr}, the specialization of $\OO_q(G)/J_{u,v}$ at $q=1$
is the coordinate ring
of $\overline {G^{u,v}}$.
Thus, we will denote $\OO_q(G)/J_{u,v}$ by $\OO_q(\overline {G^{u,v}})$
and refer to it as a \emph{quantum closed double Bruhat cell}.

In order to define the ``non-closed" quantum double Bruhat cells, we introduce
the quantum analogs of generalized minors from \cite{fz-double}.
Fix a dominant weight $\lambda \in P^+$, a pair $(u,v) \in W\times W$,
a reduced word $(i_1, \dots, i_{\ell(u)})$ for~$u$, and
a reduced word $(j_1, \dots, j_{\ell(v)})$ for~$v$.
For $k \in [1, \ell(u)]$ (resp. $k \in [1, \ell(v)]$), we define
the coroot $\eta_k^\vee$ (resp. $\zeta_k^\vee$) by setting
$\eta_k^\vee = s_{i_{\ell(u)}} \cdots s_{i_{k+1}}
\alpha_{i_k}^\vee$ (resp. $\zeta_k^\vee = s_{j_{\ell(v)}} \cdots s_{j_{k+1}}
\alpha_{j_k}^\vee$).
It is well-known that the coroots $\eta_1^\vee, \dots, \eta_{\ell(u)}^\vee$
(resp. $\zeta_1^\vee, \dots, \zeta_{\ell(v)}^\vee$) are positive
and distinct; in particular, we have $\lambda(\eta_k^\vee) \geq 0$
and $\lambda(\zeta_k^\vee) \geq 0$.
Then we define an element
$\Delta_{u \lambda, v \lambda}\in \mathcal{E}_\lambda \subset \OO_q(G)$ by
\begin{equation}
\label{eq:quantum-minor}
\Delta_{u \lambda, v \lambda} =
(F_{j_1}^{[\lambda(\zeta_1^\vee);j_1]} \cdots
F_{j_{\ell(v)}}^{[\lambda(\zeta_{\ell(v)}^\vee);j_{\ell(v)}]})
\bullet \Delta^\lambda \bullet
(E_{i_{\ell(u)}}^{[\lambda(\eta_{\ell(u)}^\vee);i_{\ell(u)}]} \cdots
E_{i_1}^{[\lambda(\eta_1^\vee);i_1]})
\end{equation}
(see \eqref{eq:divided-power}); in view of the
quantum Verma relations \cite[Proposition~39.3.7]{Lusztig}
the element $\Delta_{u \lambda, v \lambda}$ indeed depends only on
the weights $u\lambda$ and $v\lambda$, not on the choices of
$u$, $v$ and their reduced words.
It is also immediate that each quantum minor $\Delta_{\gamma,\delta}$
belongs to the graded component $\OO_q(G)_{\gamma, \delta}$, and
that it spans the one-dimensional weight space
$\mathcal{E}_\lambda \cap \OO_q(G)_{\gamma, \delta}$.
This implies that
\begin{align}
\label{eq:e-Delta}
E_i\bullet \Delta_{\gamma,\delta} =0 & \,\, \text{if} \,\, (\alpha_i\ |\ \delta)\geq 0,\\
\nonumber
F_i\bullet \Delta_{\gamma,\delta} =0 & \,\, \text{if} \,\, (\alpha_i\ |\ \delta)\leq 0;
\end{align}
\begin{align}
\label{eq:Delta-f}
\Delta_{\gamma,\delta} \bullet F_i = 0 & \,\, \text{if} \,\, (\alpha_i\ |\ \gamma)\geq 0,\\
\nonumber
\Delta_{\gamma,\delta} \bullet E_i =0 &\,\, \text{if} \,\, (\alpha_i\ |\ \gamma)\leq 0;
\end{align}

The generalized minors have the following multiplicative property:
\begin{equation}
\label{eq:minors-mult}
\Delta_{u\lambda,v\lambda}\Delta_{u\mu,v\mu} =
\Delta_{u(\lambda+\mu),v(\lambda+\mu)} \quad
(\lambda, \mu \in P^+, \,\, u, v \in W) \ .
\end{equation}
For $u=v=1$, this follows at once from \eqref{eq:highest-vector};
for general~$u$ and~$v$, \eqref{eq:minors-mult} follows
by a repeated application of the following useful lemma which is
proved by a direct calculation using \eqref{eq:coproduct} and
\eqref{eq:translation-product}.

\begin{lemma}
\label{lem:max-powers-1}
Let $f \in \OO_q(G)_{\gamma,\delta}$ and
$g \in \OO_q(G)_{\gamma',\delta'}$.
For a given~$i \in [1,r]$, suppose that~$a = \delta(\alpha_i^\vee)$
(resp.~$b = \delta'(\alpha_i^\vee)$) is
the maximal nonnegative integer such that $F_i^a \bullet f \neq 0$
(resp. $F_i^b \bullet g \neq 0$).
Then
\begin{equation}
\label{eq:max-powers-1}
(F_i^{[a;i]} \bullet f) \cdot (F_i^{[b;i]} \bullet g) =
 F_i^{[a+b;i]} \bullet (fg) \ .
\end{equation}
Similarly, if~$c = \gamma(\alpha_i^\vee)$ (resp.~$d = \gamma'(\alpha_i^\vee)$) is
the maximal nonnegative integer such that $f \bullet E_i^c  \neq 0$
(resp. $g \bullet E_i^d \neq 0$), then
\begin{equation}
\label{eq:max-powers-2}
(f \bullet E_i^{[c;i]}) \cdot (g \bullet E_i^{[d;i]}) =
(fg) \bullet E_i^{[c+d;i]} \ .
\end{equation}
\end{lemma}

The following fact can be deduced from the proof of Proposition~II.4.2 in
\cite{brown-goodearl}.

\begin{proposition}
\label{pr:modI q-commute}
For any dominant weight $\lambda\in P^+$, a pair of
Weyl group elements $u, v \in W$, and a homogeneous element
$f \in \OO_q(G)_{\gamma,\delta}$, we have
\begin{equation}
\label{eq:modI q-commute1}
f \cdot \Delta_{\lambda,v^{-1}\lambda} - q^{(\gamma\,|\,\lambda)-(\delta\,|\,v^{-1}\lambda)}
\Delta_{\lambda,v^{-1}\lambda} \cdot f \in J_{u,v} \ ,
\end{equation}
\begin{equation}
\label{eq:modI q-commute2}
\Delta_{u\lambda,\lambda} \cdot f - q^{(\gamma\,|\,u\lambda)-(\delta\,|\,\lambda)}
f \cdot \Delta_{u\lambda,\lambda} \in J_{u,v} \ .
\end{equation}
\end{proposition}

Let $\pi_{u,v}$ denote the projection
$\OO_q(G) \to \OO_q(\overline {G^{u,v}})$.
It is not hard to check that $\pi_{u,v} (\Delta_{u\lambda,\lambda})\ne 0$ and
$\pi_{u,v}(\Delta_{\lambda,v^{-1}\lambda})\ne 0$.
We can rewrite \eqref{eq:modI q-commute1} and \eqref{eq:modI q-commute2} as
\begin{align}
\label{eq:modI q-commute11}
f \cdot \pi_{u,v}(\Delta_{\lambda,v^{-1}\lambda}) &
=q^{(\gamma\,|\,\lambda)-(\delta\,|\,v^{-1}\lambda)}
\pi_{u,v}(\Delta_{\lambda,v^{-1}\lambda}) \cdot f, \\
\label{eq:modI q-commute21}
\pi_{u,v}(\Delta_{u\lambda,\lambda}) \cdot f  &
= q^{(\gamma\,|\,u\lambda)-(\delta\,|\,\lambda)}
f  \cdot \pi_{u,v}(\Delta_{u\lambda,\lambda})
\end{align}
(for $f \in \OO_q(\overline {G^{u,v}})_{\gamma,\delta}$).

In view of \eqref{eq:modI q-commute11}-\eqref{eq:modI q-commute21} and \eqref{eq:minors-mult},
for each $u,v\in W$ the set
$$D_{u,v}:=\{q^k \pi_{u,v}(\Delta_{u\lambda,\lambda})\cdot
\pi_{u,v}(\Delta_{\mu,v^{-1}\mu}): k\in \ZZ, \lambda, \mu\in P^+\}$$
is an Ore set in the Ore domain $\OO_q(\overline {G^{u,v}})$ (see Section~\ref{sec:appendix}).
This motivates the following definition.

\begin{definition}
\label{def:quantum-double-cell}
The {\it quantum double Bruhat cell}
$\OO_q(G^{u,v})$ is the localization of $\OO_q(\overline {G^{u,v}})$
by the Ore set $D_{u,v}$, that is,
$\OO_q(G^{u,v}) = \OO_q(\overline {G^{u,v}})[D_{u,v}^{-1}]$.
\end{definition}

Definition~\ref{def:quantum-double-cell} is easily seen to
coincide with the definition in \cite[Section~II.4.4]{brown-goodearl}.

\section{Cluster algebra setup in quantum double Bruhat cells}
\label{sec:quantum double Bruhat cells}

\subsection{Clusters associated with double reduced words}
\label{sec:quantum-ansatz}
Fix a pair $(u,v)\in W\times W$, and let $m = r + \ell(u) + \ell(v) = \dim G^{u,v}$.
Let $\ii=(i_1,\ldots,i_m)$ be a double word such that $(i_{r+1},\ldots, i_m)$ is a
reduced word for $(u,v)$, and $(i_1, \ldots, i_r)$ is a permutation of $[1,r]$.
For $k=1,\ldots,m$, we define the weights $\gamma_k, \delta_k \in P$ as
follows:
$$\gamma_k =
s_{-i_1} \cdots s_{-i_k}\omega_{|i_k|}, \quad
\delta_k =
s_{i_m} \cdots s_{i_{k+1}} \omega_{|i_k|}$$
(with our usual convention that $s_{-i}=1$
for $i \in [1,r]$).
Let $\Delta_{\gamma_k,\delta_k} \in \OO_q(G)$ be the corresponding
quantum minor.
Note that
$$\{\Delta_{\gamma_1,\delta_1}, \ldots, \Delta_{\gamma_r,\delta_r}\}
= \{\Delta_{\omega_1, v^{-1} \omega_1}, \ldots,
\Delta_{\omega_r, v^{-1} \omega_r}\},$$
and $\Delta_{\gamma_k,\delta_k} = \Delta_{u \omega_{|i_k|},\omega_{|i_k|}}$
whenever $k^+ = m+1$ (see Section~\ref{sec:double-words}); thus,
the only minors $\Delta_{\gamma_k,\delta_k}$ that depend on the
choice of~$\ii$ are those for which~$k$ is $\ii$-exchangeable.

\begin{theorem}
\label{th:correspondence delta}
The quantum minors $\Delta_{\gamma_k,\delta_k}$ pairwise
quasi-commute in $\OO_q(G)$.
More precisely, for $1 \leq \ell < k \leq m$, we have
\begin{equation}
\label{eq:q-commute minors}
\Delta_{\gamma_k,\delta_k} \cdot \Delta_{\gamma_\ell,\delta_\ell} =
q^{(\gamma_k\,|\,\gamma_\ell)-(\delta_k\,|\,\delta_\ell)}
\Delta_{\gamma_\ell,\delta_\ell} \cdot \Delta_{\gamma_k,\delta_k} \ .
\end{equation}
\end{theorem}

\begin{proof}
The identity \eqref{eq:q-commute minors} is a special case of the
following identity:
\begin{equation}
\label{eq:flag qcommute}
\Delta_{s's\lambda,t'\lambda} \cdot \Delta_{s'\mu,t't\mu}  =
q^{(s\lambda\,|\,\mu)-(\lambda\,|\,t\mu)}
\Delta_{s'\mu,t't\mu} \cdot \Delta_{s's\lambda,t'\lambda}
\end{equation}
for any $\lambda,\mu\in P^+$, and $s,s',t,t'\in W$ such that
$$\ell(s's)=\ell(s')+\ell(s), \quad \ell(t't)=\ell(t')+\ell(t) \ .$$
Indeed, \eqref{eq:q-commute minors} is obtained from \eqref{eq:flag qcommute}
by setting
$$\lambda = \omega_{|i_k|}, \quad \mu = \omega_{|i_\ell|},
\quad s' = s_{-i_1} \cdots s_{-i_\ell}, \quad s = s_{-i_{\ell+1}} \cdots s_{-i_k},$$
$$t' = s_{i_m} \cdots s_{i_{\max(k,r)+1}}, \quad
t =
\begin{cases}
s_{i_k} \cdots s_{i_{\max(\ell,r)+1}} & \text{if $r < k$}, \\
1 & \text{otherwise}.
  \end{cases}
$$

To prove \eqref{eq:flag qcommute}, we first consider its special
case with $s' = t' = 1$:
\begin{equation}
\label{eq:flag qcommute simple}
\Delta_{s\lambda, \lambda} \cdot \Delta_{\mu, t\mu} =
q^{(s\lambda\,|\,\mu)-(\lambda\,|\,t\mu)}
\Delta_{\mu, t\mu} \cdot \Delta_{s\lambda, \lambda}
\end{equation}
for any $\lambda,\mu\in P^+$ and $s,t \in W$.
In view of \eqref{eq:e-Delta} and \eqref{eq:Delta-f},
the minors in \eqref{eq:flag qcommute simple} satisfy
$$E_i \bullet \Delta_{s\lambda, \lambda} =
\Delta_{\mu, t\mu} \bullet F_i = 0
\quad (i \in [1,r]) \ ,$$
or equivalently,
$$E \bullet \Delta_{s\lambda, \lambda} = \varepsilon (E) \Delta_{s\lambda, \lambda}
\quad (E \in U^+), \quad
\Delta_{\mu, t\mu} \bullet F = \varepsilon (F) \Delta_{\mu, t\mu}
\quad (F \in U^-) \ .$$
Thus, \eqref{eq:flag qcommute simple} is a consequence of the following lemma.

\begin{lemma}
\label{lem:qcommute-simple-general}
Suppose the elements $f \in \OO_q(G)_{\gamma,\delta}$ and
$g \in \OO_q(G)_{\gamma',\delta'}$ satisfy
$$E \bullet f = \varepsilon (E) f \quad (E \in U^+), \quad
g \bullet F = \varepsilon (F) g \quad (F \in U^-) \ .$$
Then
\begin{equation}
\label{eq:Gauss q-commute}
fg = q^{(\gamma\,|\,\gamma')-(\delta\,|\,\delta')}gf \ .
\end{equation}
\end{lemma}

\begin{proof}
It suffices to show that both sides of \eqref{eq:Gauss q-commute}
take the same value at each element $F K_\lambda E \in U$, where
$F$ (resp.~$E$) is some monomial in $F_1, \dots, F_r$ (resp.
$E_1, \dots, E_r$).
Using \eqref{eq:translation-product} together with
\eqref{eq:coproduct}--\eqref{eq:counit} and \eqref{eq:PP-grading}, we obtain
$$(fg)(F K_\lambda E) = (E \bullet fg \bullet F)(K_\lambda) =
\sum (E_1 \bullet f \bullet F_1)(K_\lambda) \cdot
(E_2 \bullet g \bullet F_2)(K_\lambda)$$
$$= (K_{{\rm deg} E} \bullet f \bullet F) (K_\lambda) \cdot
(E \bullet g \bullet K_{{\rm deg} F}) (K_\lambda)
= q^{({\rm deg} E | \delta) + ({\rm deg} F | \gamma')}
f(F K_\lambda) \cdot g(K_\lambda E) \ ;$$
similarly,
$$(gf)(F K_\lambda E) = f(F K_\lambda) \cdot g(K_\lambda E) \ .$$
In view of \eqref{eq:UO-pairing}, $f(F K_\lambda) \neq 0$
(resp. $g(K_\lambda E) \neq 0$) implies that
${\rm deg} F = \gamma - \delta$ (resp.
${\rm deg} E = \gamma' - \delta'$).
We conclude that
$$fg = q^{(\gamma' - \delta'\, |\, \delta) + (\gamma - \delta\, |\, \gamma')} gf
=  q^{(\gamma\,|\,\gamma')-(\delta\,|\,\delta')}gf,$$
as  claimed.
\end{proof}

To finish the proof of Theorem~\ref{th:correspondence delta},
it remains to deduce \eqref{eq:flag qcommute} from \eqref{eq:flag qcommute simple}.
Remembering the definition \eqref{eq:quantum-minor}, we see that
this implication is obtained by a repeated application of the following
lemma, which is immediate from Lemma~\ref{lem:max-powers-1}.

\begin{lemma}
\label{lem:max-powers}
In the situation of Lemma~\ref{lem:max-powers-1},
suppose the elements~$f$ and~$g$ quasi-commute, i.e.,
$fg = q^k gf$ for some integer~$k$.
Then
\begin{eqnarray}
\label{eq:max-powers-qcommute}
(F_i^{[a;i]} \bullet f) \cdot (F_i^{[b;i]} \bullet g) = q^k
(F_i^{[b;i]} \bullet g) \cdot (F_i^{[a;i]} \bullet f) \ ;\\
\label{eq:max-powers-qcommute-2}
(f \bullet E_i^{[c;i]}) \cdot (g \bullet E_i^{[d;i]}) = q^k
(g \bullet E_i^{[d;i]}) \cdot (f \bullet E_i^{[c;i]}) \ .
\end{eqnarray}
\end{lemma}

This completes the proof of Theorem~\ref{th:correspondence delta}.
\end{proof}

\begin{remark}
Under the specialization $q = 1$, Theorem~\ref{th:correspondence delta}
evaluates the standard Poisson-Lie brackets between the ordinary
generalized minors.
This answer agrees with the one given in
\cite[Theorem~2.6]{kogzel}, in view of \cite[Theorem~3.1]{gsv};
in fact, Theorem~\ref{th:correspondence delta} allows one to
deduce each of these two results from another one
(see \cite[Remark~2.8]{kogzel}).
(Unfortunately, the Poisson bracket used in \cite{kogzel} and
borrowed from \cite{ks} is the opposite of the one in \cite{brown-goodearl}.)
\end{remark}

\subsection{The dual Lusztig bar-involution}
\label{sec:dual-bar-involution}
Following G.~Lusztig, we denote by $u \mapsto \overline u$
the involutive ring automorphism of $U$ such that
$$\overline q=q^{-1}, \quad \overline {E_i}=E_i, \quad
\overline {F_i}=F_i, \quad \overline {K_\mu}=K_{-\mu} \ .$$
Clearly, this involution preserves the grading
\eqref{eq:U-grading}.
Define the \emph{dual bar-involution} $f \mapsto \overline f$ on
$\OO_q(G)$ by
\begin{equation}
\label{eq:dual bar}
\overline f(u)=\overline {f(\overline u)}  \quad (u \in U) \ .
\end{equation}
This is an involutive automorphism of $\OO_q(G)$ as a
$\QQ$-vector space, satisfying
$\overline {Q f} = \overline Q \ \overline f$ for
$Q \in \QQ(q)$, where $\overline Q (q) = Q(q^{-1})$.
The definitions imply at once that
\begin{equation}
\label{eq:bar-translation}
\overline{Y \bullet f \bullet X}=
\overline Y \bullet \overline f \bullet \overline X
\quad (X,Y \in U, \, f \in \OO_q(G)) \ .
\end{equation}
It follows that
$$\overline {\OO_q(G)_{\gamma, \delta}} = \OO_q(G)_{\gamma, \delta}$$
for any $\gamma,\delta\in P$.

The dual bar-involution has the following useful multiplicative property.

\begin{proposition} For any $f \in \OO_q(G)_{\gamma,\delta}$ and
$g\in \OO_q(G)_{\gamma',\delta'}$, we have
\begin{equation}
\label{eq:bar product}
\overline {f \cdot g}=q^{(\delta\,|\,\delta')-(\gamma\,|\,\gamma')}\overline g
\cdot \overline f \ .
\end{equation}
\end{proposition}

\begin{proof}
We start with some preparation concerning ``twisted"
comultiplications in~$U$.
For a ring homomorphism $D: U \to U \otimes U$ and a ring
automorphism $\varphi$ of $U$, we define the twisted ring
homomorphism ${}^\varphi \!D: U \to U \otimes U$ by
\begin{equation}
\label{eq:twisted-D}
{}^\varphi \!D = (\varphi \otimes \varphi) \circ D \circ
\varphi^{-1} \ .
\end{equation}

In particular, we have a well defined ring homomorphism
${}^-\!\Delta: U \to U \otimes U$ corresponding to $D = \Delta$ and
$\varphi (u) = \overline u$.
Clearly, ${}^-\!\Delta$ is $\QQ(q)$-linear.

Let $\sigma:U\to U$ denote a $\QQ(q)$-linear automorphism of $U$
given by
$$\sigma(u)=q^{\frac{(\alpha\,|\,\alpha)}{2}}u K_\alpha$$
for $u\in U_\alpha$ (an easy check shows that~$\sigma$ is a ring automorphism of~$U$).
As an easy consequence of \eqref{eq:UO-pairing}, we see that
\begin{equation}
\label{eq:f-sigma}
f \circ \sigma = q^{\frac{(\gamma\,|\,\gamma)-(\delta\,|\,\delta)}{2}} f
\end{equation}
for any $f \in \OO_q(G)_{\gamma,\delta}$.

Let ${}^\sigma \!\Delta^{\rm op}:U\to U\otimes U$ be the
$\QQ(q)$-algebra homomorphism defined as in \eqref{eq:twisted-D}
with $\varphi = \sigma$ and $D = \Delta^{\rm op}$,
the \emph{opposite comultiplication} given by
$\Delta^{\rm op}= P \circ \Delta$, where
$P(X \otimes Y) = Y \otimes X$.
We claim that
\begin{equation}
\label{eq:bar op}
{}^-\!\Delta = {}^\sigma \!\Delta^{\rm op} \ ;
\end{equation}
indeed, both sides are $\QQ(q)$-algebra homomorphisms
$U\to U\otimes U$, so it suffices to show that they take the same
value at each of the generators $E_i$, $F_i$, and $K_\lambda$,
which is done by a straightforward calculation.

Now everything is ready for the proof of \eqref{eq:bar product}, which we
prefer to prove in an equivalent form:
$\overline {\overline f \cdot \overline g} =
q^{(\delta\,|\,\delta')-(\gamma\,|\,\gamma')}gf$.
Indeed, combining the definitions with \eqref{eq:bar op}
and \eqref{eq:f-sigma}, we obtain:
$$\overline {\overline f \cdot \overline g}(u) =
(f \otimes g)({}^-\!\Delta (u)) =
(f \otimes g)({}^\sigma \!\Delta^{\rm op} (u)) =
(((g \circ \sigma)\cdot (f \circ \sigma)) \circ \sigma^{-1})(u)$$
$$= q^{\frac{(\gamma\,|\,\gamma)-(\delta\,|\,\delta) +
(\gamma'\,|\,\gamma')-(\delta'\,|\,\delta') -
(\gamma+\gamma'\,|\,\gamma+\gamma')+
(\delta+\delta'\,|\,\delta+\delta')}{2}} (gf)(u) =
q^{(\delta\,|\,\delta') - (\gamma\,|\,\gamma')} (gf)(u) \ ,$$
as claimed.
\end{proof}


\begin{proposition}
\label{le:bar invariant minor}
Every quantum minor $\Delta_{\gamma,\delta}$ is invariant under
the dual bar-involution.
\end{proposition}

\begin{proof} First, we note that $\overline {\Delta^\lambda}=\Delta^\lambda$:
this is a direct consequence of \eqref{eq:highest-vector}.
The general statement
$\overline {\Delta_{\gamma,\delta}}=\Delta_{\gamma,\delta}$
follows from \eqref{eq:quantum-minor}
together with \eqref{eq:bar-translation} and the observation that
all divided powers of the elements $E_i$ and $F_i$ in $U$ are invariant under
the Lusztig involution.
\end{proof}

Let $\ii$ and the corresponding quantum minors
$\Delta_{\gamma_k,\delta_k}$ for $k=1,\ldots,m$ be as in
Section~\ref{sec:quantum-ansatz}.
Generalizing Proposition~\ref{le:bar invariant minor}, we now
prove the following.

\begin{proposition}
\label{pr:bar-invariance-monomial}
Every monomial
$\Delta_{\gamma_1,\delta_1}^{a_1}\cdots \Delta_{\gamma_m,\delta_m}^{a_m}$
is invariant under the dual bar-involution.
\end{proposition}

\begin{proof}
Using Propositions~\ref{eq:bar product}, \ref{le:bar invariant minor},
and Theorem~\ref{th:correspondence delta}, we obtain
$$\overline {\Delta_{\gamma_1,\delta_1}^{a_1}\cdots
\Delta_{\gamma_m,\delta_m}^{a_m}}=
q^{\sum_{\ell < k} a_k a_\ell
((\delta_k\,|\,\delta_\ell)-(\gamma_k\,|\,\gamma_\ell))}
\Delta_{\gamma_m,\delta_m}^{a_m}\cdots
\Delta_{\gamma_1,\delta_1}^{a_1} = \Delta_{\gamma_1,\delta_1}^{a_1}\cdots
\Delta_{\gamma_m,\delta_m}^{a_m} \ ,$$
as claimed.
\end{proof}

Note that the projection $\pi_{u,v}: \OO_q(G) \to \OO_q(\overline {G^{u,v}})$
gives rise to a well-defined dual bar-involution on $\OO_q(\overline {G^{u,v}})$
given by $\overline {\pi_{u,v}(f)} = \pi_{u,v}(\overline f)$
(indeed, the Lusztig involution preserves $U_{u,v}$ so its dual
preserves $J_{u,v} = \ker \ \pi_{u,v}$).

\begin{proposition}
\label{pr:cluster-monomials-Guv}
The monomials $\pi_{u,v}(\Delta_{\gamma_1,\delta_1})^{a_1}\cdots
\pi_{u,v}(\Delta_{\gamma_m,\delta_m})^{a_m}$ are linearly
independent over $\QQ(q)$, and each of them is invariant
under the dual bar-involution in $\OO_q(\overline {G^{u,v}})$.
\end{proposition}

\begin{proof}
The linear independence is clear because it holds under the
specialization $q = 1$.
The invariance under the dual bar-involution is immediate from
Proposition~\ref{pr:bar-invariance-monomial}.
\end{proof}

\subsection{Connections with cluster algebras}
\label{sec:q cluster setup for double cells}
As in Section~\ref{sec:quantum-ansatz},
let $\ii=(i_1,\ldots,i_m)$ be a double word such that $(i_{r+1},\ldots, i_m)$ is a
reduced word for $(u,v)$, and $(i_1, \ldots, i_r)$ is a permutation of $[1,r]$.
Let $\Lambda(\ii)$ (resp.~$\Sigma(\ii)$) be the skew-symmetric (resp. symmetric) integer
$m \times m$ matrix  defined by \eqref{eq:lambda-sigma-ii}.
We identify $\Lambda(\ii)$ with the corresponding skew-symmetric bilinear form
on $L = \ZZ^m$, and consider the based quantum torus $\TT(\Lambda(\ii))$
associated with $L$ and $\Lambda(\ii)$ according to
Definition~\ref{def:quantum-torus}.
For $k = 1, \dots, m$, we denote $X_k = X^{e_k}$, where
$\{e_1, \dots, e_m\}$ is the standard basis in $\ZZ^m$.
Let $\FF$ be the skew-field of fractions of $\TT(\Lambda(\ii))$,
and let $M : \ZZ^m \to \FF - \{0\}$ be the toric frame such that
$M(e_k) = X_k$ for $k \in [1,m]$ (see Definition~\ref{def:toric-frame}
and Lemma~\ref{lem:frame-positive-part}).

On the other hand, let $\OO_{q^{1/2}}(G^{u,v})$ denote the algebra
obtained from $\OO_q(G^{u,v})$ by extending the scalars from
$\QQ(q)$ to $\QQ(q^{1/2})$.
Let $\TT_\ii \subset \OO_{q^{1/2}}(G^{u,v})$ denote the quantum
subtorus of $\OO_{q^{1/2}}(G^{u,v})$ generated by the elements
$\pi_{u,v}(\Delta_{\gamma_1,\delta_1}), \dots,
\pi_{u,v}(\Delta_{\gamma_m,\delta_m})$ (see
Proposition~\ref{pr:cluster-monomials-Guv}).

\begin{proposition}
\label{pr:torus-identification}
\begin{enumerate}
\item
The correspondence
$X_k \mapsto \pi_{u,v}(\Delta_{\gamma_k,\delta_k})
\,\, (k \in [1,m])$
extends uniquely to a $\QQ(q^{1/2})$-algebra
isomorphism $\varphi: \TT(\Lambda(\ii)) \to \TT_\ii$.
\item
The isomorphism $\varphi$ transforms
the twisted bar-involution
$X \mapsto \overline {X}^{(\Sigma(\ii))}$ on $\TT(\Lambda(\ii))$
(see \eqref{eq:twisted-bar-involution}) into the dual
bar-involution on $\TT_\ii$ (see
Section~\ref{sec:dual-bar-involution}).
\end{enumerate}
\end{proposition}

\begin{proof}
(1) Comparing \eqref{eq:Xi-q-com} with \eqref{eq:q-commute minors},
and using Proposition~\ref{pr:cluster-monomials-Guv}, we see that
it suffices to prove the following:
\begin{equation}
\label{eq:lambda-gamma-delta}
\lambda_{k\ell}(\ii) = (\gamma_k\,|\,\gamma_\ell)-(\delta_k\,|\,\delta_\ell)
\end{equation}
for $1 \leq \ell < k \leq m$.
Remembering \eqref{eq:lambda-sigma-ii} and \eqref{eq:eta-kl}, we
obtain (for $\ell < k$):
$$(\gamma_k\,|\,\gamma_\ell)-(\delta_k\,|\,\delta_\ell) =
(s_{-i_1} \cdots s_{-i_k}\omega_{|i_k|}\,|\,
s_{-i_1} \cdots s_{-i_\ell}\omega_{|i_\ell|})$$
$$- (s_{i_m} \cdots s_{i_{k+1}} \omega_{|i_k|}\,|
\,s_{i_m} \cdots s_{i_{\ell+1}} \omega_{|i_\ell|})$$
$$= (s_{-i_{\ell+1}} \cdots s_{-i_k}\omega_{|i_k|}\,|\,
\omega_{|i_\ell|}) -
(\omega_{|i_k|}\,|
\,s_{i_k} \cdots s_{i_{\ell+1}} \omega_{|i_\ell|})$$
$$= (\pi_- [\ell^+, k] \omega_{|i_k|} - \pi_+ [\ell^+, k] \omega_{|i_k|}
 \ | \ \omega_{|i_\ell|}) = \eta_{k \ell^+} =
\lambda_{k \ell}(\ii) \ ,$$
as required.

(2) This is a direct consequence of
\eqref{eq:twisted-bar-involution}, \eqref{eq:M-Xi}
and the last statement in
Proposition~\ref{pr:cluster-monomials-Guv}.
\end{proof}

In view of Proposition~\ref{pr:torus-identification}, the
isomorphism $\varphi: \TT(\Lambda(\ii)) \to \TT_\ii$ extends
uniquely to an injective homomorphism of
skew-fields of fractions $\FF \to \FF(\OO_{q^{1/2}}(G^{u,v}))$,
which we will denote by the same symbol~$\varphi$.
Let $\mathcal{U}(M, \tilde B(\ii)) \subset \FF$ be the upper
cluster algebra associated according to \eqref{eq:upper-bound}
with the toric frame~$M$ and the matrix $\tilde B(\ii)$ given by
\eqref{eq:tildeB-entries}.
We can now state the following conjecture whose classical
counterpart is \cite[Theorem~2.10]{bfz-clust3}.

\begin{conjecture}
\label{con:birational isomorphism}
The homomorphism
$\varphi:\FF \to \FF(\OO_{q^{1/2}}(G^{u,v}))$
is an isomorphism of skew fields; furthermore, it
restricts to an isomorphism of $\QQ(q^{1/2})$-algebras
$\mathcal{U}(M, \tilde B(\ii)) \to \OO_{q^{1/2}}(G^{u,v})$.
\end{conjecture}

For instance, if $G = SL_3$, and $G^{u,v}$ is the open double
Bruhat cell in~$G$ (i.e., $u = v = w_0$) then we conjecture
that $\OO_{q^{1/2}}(G^{u,v})$ identifies with the quantum upper
cluster algebra associated with the compatible pair
$(\Lambda, \tilde B)$ in Examples~\ref{ex:sl3-quantized-matrices}
and \ref{ex:ii-standard-sl3}.


\section{Appendix: Ore domains and skew fields of fractions}
\label{sec:appendix}

Let $R$ be a \emph{domain}, i.e., an associative ring with unit
having no zero-divisors.
As in \cite[A.2]{joseph}, we say that $R$ is an Ore domain
if is satisfies the (left) Ore
condition: $aR \cap bR \neq \{0\}$ for any non-zero $a, b \in R$.
Let $\FF(R)$ denote the set of \emph{right fractions} $ab^{-1}$ with
$a,b\in R$, and $b\neq 0$; two such fractions $ab^{-1}$ and $cd^{-1}$
are identified if $af=cg$ and $bf=dg$ for some non-zero $f,g\in R$.
The ring~$R$ is embedded into $\FF(R)$ via $a \mapsto a \cdot 1^{-1}$.
It is well known that if~$R$ is an Ore domain then the addition
and multiplication in~$R$ extend to $\FF(R)$ so that
$\FF(R)$ becomes a skew-field.
(Indeed, we can define
$$a b^{-1} + c d^{-1} = (ae + cf) g^{-1},$$
where non-zero elements $e, f$, and $g$ of $R$ are chosen so that
$be = df = g$; similarly,
$$a b^{-1} \cdot c d^{-1} = ae \cdot (df)^{-1},$$
where non-zero $e, f \in R$ are chosen so that $cf = be$.)

A subset $D \subset R - \{0\}$ is called an Ore set if~$D$ is a
multiplicative monoid with unit satisfying~$dR = Rd$ for all~$d \in D$.
One checks easily that if~$D$ is an Ore set, then the set of right
fractions $R[D^{-1}] = \{ad^{-1}: a \in R, \, d \in D\}$ is a
subring of $\FF(R)$, called the localization of~$R$ by~$D$.

We now present a helpful sufficient condition for a domain
to be an Ore domain.
Suppose that~$R$ is an algebra over a field~$k$ with an
increasing filtration $\left(k \subset R_0 \subset R_1 \subset \cdots \right)$,
where each $R_i$ is a finite dimensional $k$-vector space,
$R_i R_j \subset R_{i+j}$, and $R = \cup R_i$.
We say that~$R$ has polynomial growth if $\dim R_n \le P(n)$ for all $n \geq 0$,
where $P(x)$ is some polynomial.
The following proposition is well known (see, e.g.,
\cite{ber,malikov}); for the convenience of the reader,
we will provide a proof.

\begin{proposition}
Any domain of polynomial growth is an Ore domain.
\end{proposition}

\begin{proof}
Assume, on the contrary, that $aR \cap bR = \{0\}$ for some non-zero $a, b \in R$.
Choose~$i \geq 0$ such that $a,b\in R_i$.
Then, for every $n \geq 0$, the $k$-subspaces $aR_n$ and $bR_n$ of
$R_{i+n}$ are disjoint, hence
$$\dim R_{i+n}\geq \dim aR_n + \dim bR_n \geq 2\dim R_n \ .$$
Iterating this inequality, we see that $\dim R_{mi}\geq 2^m$ for $m \ge 0$,
which contradicts the assumption that~$R$ has polynomial growth.
\end{proof}

As a corollary, we obtain that any based quantum torus
$\TT(\Lambda)$ (see Definition~\ref{def:quantum-torus})
is an Ore domain, as well as the quotient of the quantized
coordinate ring $\OO_q(G)$ (see Section~\ref{sec:OqG}) by any
prime ideal~$J$.
Indeed, both $\TT(\Lambda)$ and $\OO_q(G)/J$ are easily seen to have
polynomial growth (e.g., for $R = \OO_q(G)/J$, take $R_n$ as the
$\QQ(q)$-linear span of all products of at most~$n$ factors, each
of which is the projection of one of the generators $E_i, F_i$, or
$K_\lambda$).

We conclude with a description of the two-sided ideals in $\TT = \TT(\Lambda)$.
The following proposition is well known to the experts; it was
shown to us by Maria Gorelik.

\begin{proposition}
\label{pr:center}
\begin{enumerate}
\item
The center $Z$ of $\TT = \TT(\Lambda)$ is a free
$\ZZ[q^{\pm 1/2}]$-module with the basis
$\{X^f: f \in \ker \Lambda\}$.
Thus, $Z$ is the Laurent polynomial ring over $\ZZ[q^{\pm 1/2}]$
in $r$ independent commuting variables, where
$r = {\rm rk}(\ker \Lambda)$.

\item  The correspondence $J \mapsto I = \TT J = J \TT$
gives a bijection between the ideals in $Z$
and the two-sided ideals in $\TT$.
The inverse map is given by $I \mapsto J = I \bigcap Z$.

\item The correspondence $J \mapsto I$ in {\rm (2)}
sends intersections to intersections.
In particular, if $z_1$ and $z_2$ are relatively
prime  in~$Z$, then $\TT z_1 \cap \TT z_2 = \TT z_1 z_2$.
\end{enumerate}
\end{proposition}

\begin{proof}
We start with a little preparation.
Let $L^* = {\rm Hom}(L, \ZZ)$ be the dual lattice.
For $\xi \in L^*$, we set
\begin{equation}
\label{eq:Lstar-grading}
\TT_\xi = \{X \in T: X^e X X^{-e} = q^{\xi(e)} X
\,\, \text{for} \,\, e \in L\} \ .
\end{equation}
This makes~$\TT$ into a $L^*$-graded algebra:
the decomposition $\TT = \oplus_{\xi \in L^*} \TT_\xi$ is clear
since, in view of \eqref{eq:multiplication-F-com},
\begin{equation}
\label{T-xi}
\text{$\TT_\xi$ is a $\ZZ[q^{\pm 1/2}]$-module with the basis
$\{X^f: \xi_f = \xi\}$,}
\end{equation}
where  $\xi_f(e) = \Lambda(e,f)$.
It follows that
\begin{equation}
\label{eq:T-xi-isom}
\text{the multiplication by $X^f$ gives an isomorphism
$\TT_\xi \to \TT_{\xi + \xi_f}$.}
\end{equation}

In view of \eqref{eq:Lstar-grading}, we have $Z = \TT_0$.
Thus, assertion (1) is a special case of \eqref{T-xi}.
To prove (2), it is enough to note that every two-sided ideal~$I$
of~$\TT$ is $L^*$-graded, and, in view of \eqref{eq:T-xi-isom},
the multiplication by any $X^f$ restricts to an isomorphism
$I \bigcap \TT_\xi \to I \bigcap \TT_{\xi + \xi_f}$.
Finally, (3) is immediate from (2): since the correspondence
$I \mapsto J = I \bigcap Z$ sends intersections to intersections,
the same is true for the inverse correspondence.
\end{proof}

\section*{Acknowledgments}

We thank Maria Gorelik for telling us about
Proposition~\ref{pr:center}.
Part of this work was done while one of the authors (A.Z.) was visiting
the Warwick Mathematics Institute in April 2004; he thanks Dmitriy Rumynin
for his kind hospitality, and Ken Brown (Glasgow) for
clarifying some issues on quantum groups.


\begin{thebibliography}{xxx}

\bibitem{ber} A.~Berenstein, Group-like elements in quantum groups and Feigin's
conjecture,  to appear in {\sl Journal of Algebra}.

\bibitem{bfz-clust3}
A.~Berenstein, S.~Fomin and A.~Zelevinsky,
Cluster algebras III: Upper bounds and double Bruhat cells,
to appear in {\sl Duke Math. J.}

\bibitem{brown-goodearl}
K.~Brown and K.~Goodearl,  \textsl{Lectures on algebraic quantum groups},
Birkh\"auser, 2002.

\bibitem{dc-pr}
C.~De Concini and C.~Procesi, Quantum Schubert cells and
representations at roots of $1$, {\sl in:}  Algebraic groups and
Lie groups, {\sl Austral. Math. Soc. Lect. Ser.} \textbf{9},
Cambridge Univ. Press, Cambridge, 1997, 127--160.

\bibitem{FG1}
V.~V.~Fock and A.~B.~Goncharov,
Moduli spaces of local systems and higher Teichmuller theory,
\texttt{math.AG/0311149}.

\bibitem{FG2}
V.~V.~Fock and A.~B.~Goncharov,
Cluster ensembles, quantization and the dilogarithm,
\texttt{math.AG/0311245}.

\bibitem{fz-double}
S.~Fomin and A.~Zelevinsky, Double Bruhat cells and total
positivity, {\sl J.~Amer.\ Math.\ Soc.} \textbf{12} (1999), 335--380.

\bibitem{fz-clust1}
S.~Fomin and A.~Zelevinsky,
Cluster algebras~I: Foundations,
\textsl{J.~Amer.\ Math.\ Soc.} \textbf{15} (2002), 497--529.

\bibitem{fz-laurent}
S.~Fomin and A.~Zelevinsky,
The Laurent phenomenon,
\textsl{Adv.\ Applied Math.} \textbf{28} (2002), 119--144.

\bibitem{fz-clust2}
S.~Fomin and A.~Zelevinsky, Cluster algebras~II: Finite type
classification, \textsl{Invent.\ Math.}, \textbf{154} (2003), 63--121.

\bibitem{gsv}
M.~Gekhtman, M.~Shapiro, and A.~Vainshtein, Cluster algebras
and Poisson geometry, \textsl{Moscow Math.~ J.}, \textbf{3}, No.3 (2003).

\bibitem{gsv2}
M.~Gekhtman, M.~Shapiro, and A.~Vainshtein,
Cluster algebras and Weil-Petersson forms, \texttt{math.QA/0309138}.

\bibitem{malikov}
K.~Iohara and F.~Malikov,
Rings of skew polynomials and Gel'fand-Kirillov conjecture for quantum
groups, \textsl{Commun. Math. Phys.} \textbf{164} (1994), 217--238.

\bibitem{joseph}
A.~Joseph, \textsl{Quantum groups and their primitive ideals},
Ergebnisse der Math. (3) 29, Springer-Verlag, Berlin, 1995.

\bibitem{kac}
V.~Kac, \textsl{Infinite dimensional Lie algebras}, 3rd edition,
Cambridge University Press, 1990.

\bibitem{kogzel}
M.~Kogan and A.~Zelevinsky,
On symplectic leaves and integrable systems in
standard complex semisimple Poisson-Lie groups,
\textsl{Intern.\ Math.\ Res.\ Notices} 2002, No.32, 1685--1702.

\bibitem{ks}
L.~Korogodski and Y.~Soibelman,
\textsl{Algebras of functions on quantum groups. Part I.}
Mathematical Surveys and Monographs, \textbf{56}, American Mathematical
Society, Providence, RI, 1998.

\bibitem{Lusztig}
G.~Lusztig, \textsl{Introduction to quantum groups}, Progress in Mathematics,
\textbf{110}, Birkh\"auser Boston, 1993.

\bibitem{Lusztig-problems}
G.~Lusztig, Problems on canonical bases, {\sl in:}  Algebraic groups and
their generalizations: quantum and infinite-dimensional methods
(University Park, PA, 1991), {\sl Proc. Sympos. Pure Math.}, \textbf{56}, Part 2,
Amer. Math. Soc., Providence, RI, 1994, 169--176.

\bibitem{sz}
P.~Sherman and A.~Zelevinsky,
Positivity and canonical bases in rank 2 cluster algebras of finite
and affine types, to appear in \textsl{Moscow Math. J.}


\bibitem{z-imrn}
A.~Zelevinsky, Connected components of real double Bruhat cells,
\textsl{Intern.\ Math.\ Res.\ Notices} 2000, No. 21, 1131--1153.

\end{thebibliography}
\end{document}